
\input amstex
\documentstyle{amsppt}
\magnification = 1200
\pagewidth {6.5 true in}
\pageheight {8.5 true in }
\hcorrection {-.2 true in}
\vcorrection {+.1 true in}


\def\sA{{\Cal A}}
\def\sB{{\Cal B}}
\def\sC{{\Cal C}}
\def\sD{{\Cal D}}
\def\sE{{\Cal E}}
\def\sF{{\Cal F}}

\def\sI{{\Cal I}}

\def\sM{{\Cal M}}

\def\sP{{\Cal P}}

\def\sV{{\Cal V}}
\def\sW{{\Cal W}}
\def\nbb{\Bbb N}

\def\Ad{{\roman{Ad}}}

\def\lcm{{\roman{lcm}}}
\def\len{{\roman{len}\thinspace}}

\def\cstar{$\roman{C}^{*}$}

\def\stbar{{\thinspace |\thinspace }}
\def\indlimit{{\displaystyle \lim_{\longrightarrow} }}
\def\ldotss{\mathinner{\ldotp\ldotp\ldotp\ldotp}}
\def\lb{\lbrack}
\def\rb{\rbrack}
\def\ra{R(\sA)}
\def\nda{N_{\sD}(\sA)}
\def\resultspace{{\vskip 0.075in}}

\def\braks#1{\lb \,#1\,\rb}
\def\chisub#1 {{\chi_{\lower2.5pt\hbox{$\scriptstyle #1$}}}}

\define\row#1#2{(#1_1,\dots,#1_{#2})}
\define\rowh#1#2{(#1_1,\dots,#1_{#2},h)}
\def\tupl#1#2{(#1_0,\ldots,#1_{#2-1})}
\def\tuplone#1#2{(1,#1_1,\ldots,#1_{#2-1})}
\def\system#1#2#3#4{T_{#1} @> \phi_1 >> T_{#2} @> \phi_2 >> T_{#3} @> \phi_3 >>
 \dots @>>> #4}

\def\symbarrow#1{\smash{\mathop{\longrightarrow}\limits^{#1}}}
\def\symbdown#1{\Big\downarrow\rlap{$\vcenter{\hbox{$\scriptstyle #1$}}$}}
\def\symbup#1{\Big\uparrow\rlap{$\vcenter{\hbox{$\scriptstyle #1$}}$}}
\def\symbupdown#1{\Big\updownarrow\rlap{$\vcenter{\hbox{$\scriptstyle #1$}}$}}

\def\symbse#1{\!\!\!\!\searrow\rlap{$\vcenter{\hbox{$\scriptstyle #1$}}$}}

\topmatter
\title Order Preservation in Limit Algebras \endtitle
\author Allan P. Donsig \\  Alan Hopenwasser \endauthor
\affil Lancaster University\\
       University of Alabama
\endaffil
\address Dept. of Mathematics, Lancaster University,
      Lancaster, England, LA1 4YF  \endaddress
\email a.donsig\@lancaster.ac.uk \endemail
\address  Dept. of Mathematics, University of Alabama
      Tuscaloosa, Alabama, U.S.A., 35487 \endaddress
\email  ahopenwa\@ua1vm.ua.edu \endemail
\thanks
The first author acknowledges partial
support from an NSERC postdoctoral fellowship;
the second author acknowledges partial
support from an Alabama EPSCoR travel grant.
\endthanks
\endtopmatter

A limit algebra is the inductive limit of a system of the form:
$$
A_1 @> \alpha_1 >> A_2 @> \alpha_2 >> A_3 @> \alpha_3 >>
A_4 @> \alpha_4 >> \cdots
$$
where each $A_i$ is a member of some  class of
finite dimensional algebras and
each $\alpha_i$ is an  injective algebra homomorphism,
possibly satisfying additional specified properties.
Limit algebras have become an important source for many varied
examples of norm-closed non-self-adjoint algebras
\cite{B,HL,HPo,MS2,MS3,PePW1,PePW2,PeW,PW,Po2,Po4,SV,T,V1}.

The most restrictive non-self-adjoint context
 is for each $A_i$ to be the algebra of
upper-triangular $n_i$ by $n_i$ matrices for some sequence $n_i$.
A more general context is for each $A_i$ to be the upper-triangular
matrices in a finite-dimensional \cstar-algebra, i.e., a
 direct sum of full matrix algebras.
The most general building blocks which have been fruitful to date
are digraph algebras (finite-dimensional CSL algebras),
as described in~\cite{Po2} or Section~{6.7} of~\cite{Po4}.

There is also a range of possible assumptions for the maps.
We will assume that each map is unital, $*$-extendible and regular.
A regular map is one that maps matrix units to sums of matrix units.
This assumption ensures that there are enough partial isometries in $\sA$
which normalize $\sA \cap \sA^*$ to span $\sA$.
For a discussion of limit algebras where the maps are not regular,
see~\cite{Po3}.

Using the $*$-extendibility of the maps, we can conclude that the limit
algebra is contained, in the first context, in a UHF \cstar-algebra
and in the second, in an AF \cstar-algebra.
In fact, the \cstar-algebra generated by such a limit algebra
is the \cstar-envelope of the limit algebra.

Much of the literature has dealt with special families of limit algebras;
these usually arise in one of two ways:
\roster
\item all limit algebras possessing some natural intrinsic property,
\item all limit algebras arising from direct systems with a particular class
      of embeddings.
\endroster
This paper will deal primarily with the second situation.
The focus is on embeddings which are order preserving,
that is, which preserve the natural ordering on the diagonal (this is
defined precisely in Section~2).
In \cite{Po4} and \cite{Po5} such embeddings are called strongly regular
but the term order preserving seems more natural to us.

Two issues immediately present themselves:
\roster
\item"$\bullet$" Find intrinsic properties which characterize the family under
      consideration,
\item"$\bullet$" Classify direct systems, i.e., when do two different
      presentations yield the same algebra
      (up to isometric isomorphism) ?
\endroster
We will consider both of these questions.
The first problem has been solved elsewhere
for nest embeddings \cite{HPe} and mixing embeddings~\cite{Do};
the second question has been answered for standard embeddings
\cite{B,PePW1,Po1}, refinement embeddings \cite{PePW1,Po1}, and
alternation embeddings~\cite{HPo,P}.
The classification results in this paper subsume those for algebras with
refinement, standard and alternation embeddings.
\par
There are a number of other classification theorems
in the literature.  Algebras based on two special classes
of nest embeddings, refinement with twist embeddings and
homogeneous embeddings based on the backshift, were
classified in \cite{HPo}.  In \cite{Po3}, Power classified
an uncountable family of algebras based on nest embeddings
which are not regular.  The resulting nests generate masas which are not
canonical.  (Indeed, they generate singular masas.)  Another
classification theorem for regular, non-$*$-extendible embeddings
between certain digraph algebras can be found in \cite{Po2}.
See also \cite{Po7} for an additional example of a classification
theorem.

In answering the first question, it is natural to consider
elements of the limit algebra which preserve this diagonal order.
These are also of interest in the study of product-type cocycles;
see Lemmas~5.5 and~5.6 of \cite{V1}.
In the more general context of subalgebras of groupoid \cstar-algebras,
the concept of a monotone $G$-set is equivalent to that of an
order-preserving element; see page~57 of \cite{MS1}.
We thank Paul Muhly for pointing this out to us.

The authors would like to thank Steve Power for
several helpful conversations on the subject matter of
this paper and Dave Larson for assistance in facilitating
their collaboration.
\par
\subhead Summary of Paper \endsubhead
In Section~1, we introduce the usual examples and notation for limit algebras;
in Section~2, we define the order preserving normalizer, order preserving
embeddings and locally order preserving embeddings.
The main results of Section~2 are characterizations of locally order
preserving embeddings between $T_n$'s (Lemma~2), order preserving
embeddings between $T_n$'s (Theorem~5), and order preserving
embeddings between direct sums of $T_n$'s (Theorem~6).

In Section~3, we concentrate on the $T_n$ context and characterize
the spectra of such triangular AF algebras where the
embeddings are locally order preserving (Theorem~7) and where compositions
of the embeddings are locally order preserving (Theorem~8).

Section~4 describes in detail the spectrum of a TAF algebra
with order preserving embeddings between $T_n$'s.
In particular, using arguments similar to those in \cite{HPo} we
obtain the first step in our classification of
algebras which are limits of order preserving systems (Theorem~13).
The spectral description is also useful in constructing explicit
cocycles for these algebras, which shows these algebras and those
with order preserving embeddings through direct sums of $T_n$'s
are analytic (Theorems~14 and~15).

In Section~5, we establish various equivalent conditions for a subset
of the normalizer to generate the algebra (Proposition~17).
Applying this to the order preserving normalizer, we show that for a
triangular AF algebra, the order preserving normalizer generates the
algebra if, and only if,
 it has a presentation $\indlimit(A_i,\alpha_i)$
where $\alpha_j \circ \alpha_{j-1} \circ \cdots \circ \alpha_i$
is locally order preserving for each $i$ and $j$.

In Section~6, we state and prove a theorem relating isometric
isomorphisms and intertwining diagrams for inductive limits;
this type of theorem has appeared implicitly in a paper of Davidson
and Power, \cite{DPo}, and in slightly different forms as
Theorem~4.6 of \cite{V1} and Corollary~1.14 of \cite{PeW}.
We use this theorem to give simple proofs of several known
classification theorems and to extend a recent result of
Poon and Wagner \cite{PW}.

To apply the intertwining diagram theorem to classifying algebras
with order preserving embeddings between $T_n$'s, we need a unique
factorization theorem for order preserving embeddings between $T_n$'s.
This is obtained in Section~7.
Finally, in Section~8 we put together the results obtained in
previous sections to give a classification theorem for
algebras with order preserving embeddings between $T_n$'s.

\head 1 Preliminaries \endhead

Recall that a \cstar-algebra is {\it approximately finite} (AF) if
there is a nested sequence of finite-dimensional \cstar-algebras
whose closed union is the original \cstar-algebra.
Given an AF \cstar-algebra $\sC$ and a maximal abelian self-adjoint
subalgebra (masa) $\sM \subseteq \sC$,
we call $\sM$ a {\it canonical masa} if
 there is a nested sequence of
finite-dimensional \cstar-subalgebras of $\sC$, $\{ C_i \}$, so
that
\roster
\item $\sC = \overline{ \cup_i C_i }$,
\item $M_i = C_i \cap \sM$ is a masa in $C_i$ for each $i$, and
\item $N_{M_i}(C_i) \subseteq N_{M_{i+1}}(C_{i+1})$ for each $i$,
\endroster
where
$$ N_X(Y) = \{ y \in Y \stbar y \text{ is a partial isometry, and }
   y^* x y, y x y^* \in X \text{ for all } x \in X \}.$$
If $D_n$ is the algebra of diagonal $n \times n$ matrices and
$T_n$ is the algebra of upper-triangular $n \times n$ matrices, then
$N_{D_n}(T_n)$ consists of upper-triangular matrices with entries
either $0$ or of absolute value $1$ such that each row or column has
at most one non-zero entry.

We can now define a {\it triangular} AF (TAF) algebra, $\sA$, to be
a norm-closed subalgebra of an AF \cstar-algebra so that $\sA \cap \sA^*$
is a canonical masa in the AF \cstar-algebra.
A {\it canonical algebra} is a norm-closed subalgebra of an AF
\cstar-algebra that contains a canonical masa.

If we let $A_i$ denote $\sA \cap C_i$ and $\alpha_i$ denote the injection
map from $C_i$ to $C_{i+1}$ restricted to $A_i$, then we have a system
where each $A_i$ is finite-dimensional
 and each $\alpha_i$ is $*$-extendible.
We will call the system
$$
A_1 @> \alpha_1 >> A_2 @> \alpha_2 >> A_3 @> \alpha_3 >>
A_4 @> \alpha_4 >> \cdots \tag{1}
$$
a {\it presentation} of the TAF algebra $\sA$.

In general, each $A_i$ is not maximal as a triangular subalgebra of $C_i$
and need not be even if $\sA$
 is maximal as a triangular AF algebra in $\sC$
(see \cite{PePW1, Example 3.25} or \cite{Po4, Proposition 10.4}).
A TAF algebra is called {\it strongly maximal} if there is a sequence of
finite-dimensional \cstar-algebras $C_i$ as above so that
$A_i$ is maximal as a triangular subalgebra of $C_i$ for each $i$.
These are precisely the TAF algebras that have presentations such as (1)
with each $A_i$ either the upper-triangular $n \times n$ matrices for some $n$
or a direct sum of such.

By an {\it embedding}, we mean an injective algebra homomorphism between
triangular subalgebras of finite-dimensional \cstar-algebras that extends
to an injective $*$-homomorphism of the \cstar-algebras
and is {\it regular} in the sense that it
 maps matrix units to sums of matrix units.
Algebras built from embeddings which are not $*$-extendible
or from embeddings which are not regular have been studied in
the literature \cite{Po2,Po3,HL}; we are incorporating regularity
and $*$-extendibility into our definition since all the embeddings
that we study in this paper satisfy these two properties.  We can
thereby avoid endless repetition of these two assumptions.
If $\alpha\: A_1\rightarrow A_2$ is a regular embedding,
then $\alpha( N_{D_1}(A_1) ) \subseteq N_{D_2}(A_2) $ where
$D_i$ is the diagonal $A_i \cap A_i^*$, for $i=1,2$.
Since an embedding is $*$-extendible, if $\alpha\: T_n
\rightarrow T_m$ is an embedding, then $n$ divides $m$;
the {\it multiplicity} of $\alpha$ is the quotient $m/n$.

Two fundamental examples of embeddings have influenced much of the theory
of TAF algebras.
They are the {\it standard embedding}, $\sigma_k \: T_n \rightarrow T_{nk}$,
given by
$$ \sigma_k(A) = \bmatrix A & & & \\ & A & & \\
                        & & \ddots & \\ & & & A \endbmatrix, $$
where the righthand side is a $k$ by $k$ block matrix,
and the {\it refinement embedding} (or the canonical nest embedding),
$ \rho_k \: T_n \rightarrow T_{nk}$, given by
$$ \rho_k((a_{ij})) = \left( a_{ij} I_k \right),$$
where $I_k$ is the $k$ by $k$ identity matrix.

\head 2 Embeddings \endhead

Consider $T_n$ and its diagonal projections, denoted by $\sP(T_n)$.
The diagonal ordering on $\sP(T_n)$ (denoted $\preceq$) is a partial order
given by
$$ e \preceq f \iff
   \text{ there exists $w\in N_{D_n}(T_n)$ with $w w^* = e$, $w^* w = f$.}
$$
Notice that two comparable projections must have equal traces and
that this ordering is a total order on the minimal diagonal projections.

Each element
 $w \in N_{D_n}(T_n)$ induces a partial homeomorphism on $\sP(T_n)$,
with domain $\{ x \in \sP(T_n) \stbar x \le w w^* \}$
and range $\{ x \in \sP(T_n) \stbar x \le w^* w \}$,
given by $x \mapsto w^* x w $.

\definition{Definition}
We say that $w$ is  {\it order preserving} if
 this map preserves the diagonal ordering
restricted to its domain and range.
Define
$$ N^{op}_{D_n}(T_n) = \{ w \in N_{D_n}(T_n)
   \stbar w \text{ is order preserving} \}.$$
\enddefinition

In other words, $w \in N^{op}_{D_n}(T_n)$
if, and only if,  $x \prec y$ implies
$w^*x w \prec w^* y w $ for {\it all\/} diagonal projections $x,y$
with $x,y \le w w^*$.
Every matrix unit in $T_n$ is trivially order preserving.
The following partial isometry is not order preserving in $T_4$:
$$ \bmatrix 0 & 0 & 0 & 1 \\ & 0 & 1 & 0 \\
              &   & 0 & 0 \\ &   &   & 0 \endbmatrix.  $$

Given $x,y \in N_{D_n}(T_n)$, $x + y \in N_{D_n}(T_n)$ if, and only if,
the initial projections of $x$ and $y$ are orthogonal and
the final projections are orthogonal.
As the matrix above shows, the same is not true if we replace
$N_{D_n}(T_n)$ with $N^{op}_{D_n}(T_n)$.

The definitions below make sense for  embeddings which
are regular; $*$-extendibility is not needed.  All the order
preserving embeddings in this paper will, however, be
$*$-extendible.

\definition{Definition}
An embedding $\phi$ is {\it locally order preserving\/}
if $\phi(e)$ is order preserving for each matrix unit $e$.
An embedding $\phi$ is {\it order preserving\/}
if $\phi( N^{op}_{D_n}(T_n) ) \subset N^{op}_{D_{nk}}(T_{nk})$.
\enddefinition

 The map $\phi \colon T_2 \rightarrow T_4$ given by
$$
\phi \left( \bmatrix a & b \\ & c \endbmatrix \right) =
\bmatrix a & 0 & 0 & b \\ & a & b & 0 \\ & & c & 0 \\ & & & c \endbmatrix
$$
is an example of an embedding which
is {\it not\/} locally order preserving.
On the other hand,
 the map $\psi \colon T_{2n} \rightarrow T_{4n}$ given by
$$
\psi \left( \bmatrix A & B \\ & C \endbmatrix \right) =
\bmatrix A & 0 & B & 0 \\ & A & 0 & B \\ & & C & 0 \\ & & & C \endbmatrix
$$
with $A,C \in T_n$ and $B \in M_n$ is locally order preserving
but not order preserving for $n > 1$.
Since matrix units are in the order preserving normalizer, all
order preserving maps are locally order preserving.
Both refinement and standard embeddings are order preserving.

Much of the above discussion extends directly to direct sums of $T_n$'s.
Suppose $T = \oplus_{i=1}^a T_{m_i}$ and $D = T \cap T^*$.
We can define a diagonal order, $\preceq$, for projections
in $T$ just as before.
Notice that two minimal diagonal projections are comparable if,
and only if, they are in the same summand.

Again, each $w \in N_D(T)$ induces a partial homeomorphism on $\sP(T)$
given by $x \mapsto w^* x w$ and we can define $N^{op}_D(T)$ just
as before.
A key fact about $N^{op}_D(T)$ is that if $a,b \in N^{op}_D(T)$
are contained in {\it different\/} summands of $T$, then
$a+b \in N^{op}_D(T)$.

Just as before, we can define locally order preserving and
order preserving for an embedding
$\phi \colon \mathop{\oplus}\limits_{i=1}^a T_{m_i} \rightarrow
   \mathop{\oplus}\limits_{j=1}^b T_{n_j}$.

In either context, it is easy to see that the composition of two order
preserving embeddings is order preserving; the same is {\it not\/} true
for locally order preserving embeddings (see, for instance, the first
example on page 26 below).
Nonetheless, we have the following lemma:

\proclaim{Lemma 1}
Let $\alpha$ and $\beta$ be two embeddings.

If $\beta \circ \alpha$ is locally order preserving, then
$\alpha$ is locally order preserving.

If $\beta \circ \alpha$ is order preserving, then
$\alpha$ is order preserving.
\endproclaim

\demo{Proof}
To prove the first statement, suppose $\alpha$ is not locally order
preserving.
Then there is a matrix unit $e$ and there are
 minimal diagonal subprojections of
$\alpha(e e^*)$, $x$ and $y$, so that $x \prec y$ and
$ \alpha(e^*)x\alpha(e) \succ \alpha(e^*)y\alpha(e) $.

Since $\beta$ is regular and $*$-extendible, it is easy to see that
$$\beta(x) \prec \beta(y) \text{ and }
   \beta(\alpha(e^*)x\alpha(e)) \succ \beta(\alpha(e^*)y\alpha(e)).$$
But, as $\beta$ is a homomorphism,
$\beta(\alpha(e^*)x\alpha(e)) =
 \beta(\alpha(e^*))\beta(x)\beta(\alpha(e))$.
Thus, we have
$$\beta(x) \prec \beta(y) \text{ and }
\beta(\alpha(e^*))\beta(x)\beta(\alpha(e)) \succ
\beta(\alpha(e^*))\beta(y)\beta(\alpha(e)). $$
This shows the map $z \mapsto \beta(\alpha(e^*))z\beta(\alpha(e))$ is
not order preserving on the two subprojections $\beta(x)$ and $\beta(y)$.
Hence $\beta(\alpha(e))$ is not order preserving and the first statement
is proved.

To prove the second statement, repeat the above argument with
$\alpha$ not order preserving and $e$ in the order preserving normalizer.
\qed
\enddemo

We turn now to characterizing locally order preserving embeddings
and order preserving embeddings between $T_n$'s.

\subhead Locally Order Preserving Embeddings \endsubhead
If $\phi\colon T_n \rightarrow T_{nk}$ is locally order preserving,
then $\phi$ is determined by its action on $D_n$, and,
in particular, by its action on the minimal diagonal projections of $D_n$.

It is helpful to let $\braks{n} = \{ 1,2,\ldots,n\}$.
Clearly, we can identify
the minimal diagonal projections of $D_n$ under the diagonal ordering
with $\braks{n}$ under the usual ordering.
In the following, we use ordered pairs of integers as indices for
minimal diagonal projections in $D_{nk}$.
For the sake of clarity, we will denote such a projection by its index
alone (i.e., we will write $(i,j)$ for $e_{(i,j),(i,j)}$).
There is a bijection between the minimal diagonal projections of $D_{nk}$ and
$\braks{n} \times \braks{k}$ as follows:
For each $i \in \braks{n}$, $\phi(e_{ii})$ is the sum of $k$ minimal
diagonal projections in $D_{nk}$;
for $\phi(e_{11})$, index these projections
in the order in which they appear in the diagonal order by
$(1,1),(1,2),\ldots,(1,k)$;
for $\phi(e_{ii})$ in general, let $(i,j)$ be the image of $(1,j)$ under
conjugation by $\phi(e_{1i})$.
With this indexing, we have
$$ \phi(e_{i,j}) = \sum_{l=1}^k e_{(i,l),(j,l)}. \tag{2} $$

Notice that $\phi(e_{1j})$ is order preserving if, and only if, the diagonal
ordering restricted to $\{ (j,l) \stbar l \in [k] \}$ induces the usual order
on $[k]$.   (An obvious identification is made here.)
In general, $\phi(e_{ij})$ is order preserving if and only if, the diagonal
ordering restricted to
$\{(i,l) \stbar l \in [k] \}$ and to $\{(j,l) \stbar l \in [k] \}$
both induce the same order on $[k]$.
In particular, if $\phi$ is locally order preserving, then
the ordering on $\braks{n} \times \braks{k}$ satisfies:
$$ \aligned
i_1 < i_2 &\Longrightarrow (i_1,j) \preceq (i_2,j)
  \qquad \text{any } j \in \braks {k} \\
j_1 < j_2 &\Longrightarrow (i,j_1) \preceq (i, j_2)
  \qquad \text{any } i \in \braks {n}
\endaligned \tag{3} $$

It is straightforward to prove that:

\proclaim{Lemma 2}
There is a bijection, given by (2), between
locally order preserving embeddings
$\phi\colon T_n \rightarrow T_{nk} $ and
orderings on $\braks{n} \times \braks{k}$ that satisfy (3).
\endproclaim

This correspondence between local
order preservation and the properties of the diagonal
order on $[n] \times [k]$
depends on the $*$-extendibility of the embedding.
For example, the compression embeddings considered in \cite{HL}
are outside this framework even though they satisfy the
conditions in the definition of order preserving,
since they are not $*$-extendible.

\subhead Order Preserving Embeddings \endsubhead
Given the correspondence in Lemma~2, it is natural
to ask under what additional conditions does the diagonal ordering on
$\braks{n} \times \braks{k}$ correspond to an order preserving map?

\proclaim{Lemma 3}
A  locally order preserving embedding
$\phi\colon T_n \rightarrow T_{nk} $ is order preserving
if, and only if,
there are no $a,b \in \braks{k}$ so that
$$ (g,a) \prec (i,b) \text{ and } (h,a) \succ (j,b) $$
where $g$, $h$, $i$ and $j$ satisfy $e_{gh} + e_{ij} \in N^{op}_{D_n}(T_n)$.
\endproclaim

One consequence of this lemma is that in checking if a map such as $\phi$
is order preserving, we need only consider those elements of the
order preserving normalizer which are the sum of {\it two} matrix units.

\demo{Proof}
To establish necessity, observe that
conjugation by $\phi(e_{gh})$ maps $(h,a)$ to $(g,a)$ and
conjugation by $\phi(e_{ij})$ maps $(j,b)$ to $(i,b)$.
If there do exist $a,b \in \braks{k}$ with the given properties, then
conjugation by $\phi(e_{gh}+e_{ij})$ is
not order preserving and so $\phi$ is not order preserving.

For sufficiency, it is enough to assume $\phi$ is not order preserving
and find $a,b,g,h,i$, and $j$ satisfying the conditions.
By assumption, there is some $x \in N^{op}_{D_n}(T_n)$ so
that $\phi(x) \not\in N^{op}_{D_{n+1}}(T_{n+1})$.
Thus there are two elements  $(h,a)$ and $(j,b)$ so that conjugation
by $\phi(x)$ reverses the diagonal ordering.
Let $(g,a)$ and $(i,b)$ be their images under conjugation.
It follows that this choice of $a,b,g,h,i$ and $j$ satisfies the conditions.
\qed \enddemo

\proclaim{Lemma 4} Let $\phi\:T_n \to T_{nk}$ be an
order preserving embedding.  Let $t \in [k]$.  The diagonal
ordering on $[n] \times [k]$ satisfies
$$
(i,t) \prec (g,h) \prec (j,t) \Longrightarrow i \le g \le j
$$
for all $i$ and $j$.
\endproclaim
\demo{Proof}Assume that $(i,t) \prec (g,h) \prec (j,t)$.
If either $g<i$ or  $g>j$ then
$$
w = e_{gg} + e_{ij} \in N_{D_n}^{op}(T_n).
$$
Observe that $\phi(w)$ carries $(j,t)$ to $(i,t)$ and
$(g,h)$ to $(g,h)$.  Since $(g,h) \prec (j,t)$ and
$(g,h) \succ (i,t)$, $\phi(w)$ is not order preserving,
a contradiction.  So we must have $i \le g \le j$.
\qed
\enddemo

Given any two order preserving embeddings,
$\alpha\: T_n \rightarrow T_{na}$ and
$\beta\: T_n \rightarrow T_{nb}$,
it is easy to check that $\alpha\oplus\beta\: T_n \rightarrow T_{n(a+b)}$
is order preserving.
Since every refinement embedding is order preserving, it follows
that every direct sum of refinement embeddings is  order preserving.

The family of embeddings which are direct sums of refinement
embeddings includes all refinement embeddings (one summand
only), all standard embeddings (each refinement embedding
has multiplicity 1) and all embeddings of the form
$\sigma \circ \rho$ (the refinement embeddings in the
direct sum all have equal multiplicity).
These are the embeddings which yield the alternation
algebras studied in \cite{HPo,P,Po6}.

The next theorem shows that the direct sums of refinement
embeddings are precisely the class of order preserving
embeddings.

\proclaim{Theorem 5}
Suppose $\phi\:T_n \to T_{nk}$ is an
 embedding.
Then $\phi$ is order preserving if, and only if,
$\phi$ is a direct sum of refinement embeddings.
\endproclaim
\remark{Remark} Theorem 5 provides a description of all
order preserving embeddings in the context in which
embeddings are regular and $*$-preserving -- the context
of this paper.  Outside this setting, compression embeddings
\cite{HL} provide examples of order preserving embeddings
which are not direct sums of refinements embeddings.
\endremark
\demo{Proof}
As we have remarked earlier, it is easy to see that
a direct sum of refinement embeddings is order preserving.
To prove the converse, consider the diagonal order induced
on $[n] \times [k]$ by $\phi$ as in the discussion in the section
on locally order preserving embeddings.
The fact that the range of $\phi$ is contained in the upper triangular
matrices implies that the first element in the order is (1,1).
Let $r_1$ be the largest integer so that $(1,1), \dots, (1,r_1)$ are the
first $r_1$ elements of the diagonal order.
Lemma~4 implies that no $(g,h)$ with
$g \ge 3$ can appear in between  $(1,t)$ and $(2,t)$, so
$(2,1), \dots, (2,r_1)$ must follow in the diagonal order.
We cannot have $(2,r_1+1)$ next, for then $(2,r_1+1)$ would
precede $(1,r_1+1)$, an impossibility.
We cannot have $(1,r_1+1)$ next, as Lemma~4 implies
that between $(2,r_1)$ and $(3,r_1)$ we can only have elements
$(g,h)$ with $g=2$ or $g=3$.
Arguing in the same way, we see that
the next $r_1$ elements are $(3,1)$ through $(3,r_1)$
and so on until $(n,1)$ through $(n,r_1)$.

Upper triangularity of the image (or the conditions for
local order preservation) guarantee that the next element
is $(1,r_1 +1)$.  Now let $r_2$ be such that the diagonal
order runs $(1,r_1+1), \dots (1,r_1+r_2), (2,r_1+1)$.
We may continue the argument as before until we finally
obtain integers $r_1, \dots,r_t$ whose sum is $k$ with
the property that the embedding $\rho_{r_1} \oplus
\dots \oplus \rho_{r_t}$ induces the same diagonal order
on $[n] \times [k]$ as $\phi$ does.
Lemma~2 now yields the theorem. \qed
\enddemo

\subhead More Order Preserving Embeddings \endsubhead
In extending Theorem~5 to characterize order preserving
embeddings between direct sums of $T_n$'s, we need the
notion of an ordered Bratelli diagram, first described in
\cite{Po5}.
These diagrams play a role in \cite{HPS} and in \cite{PW};
our definitions follow those of \cite{PW}.

\definition{Definition}
Given non-empty finite sets $V$ and $W$, an {\it ordered
diagram\/} from $V$ to $W$ is a partially ordered set $E$
and maps $r \: E \rightarrow W$ and $s \: E \rightarrow V$
such that $e,e' \in E$ are comparable if and only if
$r(e) = r(e')$.
\enddefinition

The sets $V$ and $W$ are the {\it vertices\/} of the diagram
and $E$ are the {\it edges\/}.
We extend the definition slightly to describe
order preserving maps between direct sums of $T_n$'s.

\definition{Definition}
Call $(E,r,s,f)$ an {\it ordered diagram with multiplicity}
if $(E,r,s)$ is an ordered diagram as defined above, and $f$ is
a function from $E$ to $\Bbb{N}$.
\enddefinition

We call $f(e)$ the {\it multiplicity\/} of the edge $e$.

To an ordered diagram with multiplicity, we can associate
a direct sum of refinement embeddings.
Let $(E,r,s,f)$ be an ordered diagram with multiplicity
from $V = \{ 1,2,\ldots,a\}$ to $W = \{ 1, 2, \ldots, b\}$.
Given positive integers $m_1,m_2,\ldots,m_a$,
the ordered diagram with multiplicity determines a map
 $\phi \colon \oplus_{i=1}^a T_{m_i}
\rightarrow \oplus_{j=1}^b T_{n_j}$
where $ n_j = \sum_{r(e)=j} m_{s(e)} f(e)$ for each $j \in W$.
The map $\phi$ is given by
$$ \phi\left( \bigoplus_{i=1}^a t_i \right) = \bigoplus_{j=1}^b
       \left( \bigoplus_{r(e) = j} \rho_{f(e)} (t_{s(e)}) \right).
$$
It is important to stress that the inner direct sum in the definition
of $\phi$ is ordered.
That is, the diagonal ordering of the summands $\rho_{f(e)} (t_{s(e)})$ in
$T_{n_j}$ is given by the ordering of $E$ restricted to
$\{e \in E \stbar r(e)=j\}$.
This association is a slight generalization of the association between
ordered diagrams and embeddings given in \cite{Po5} and \cite{PW}.

The following definition formalizes the notion of when
two ordered diagrams with multiplicity should be
considered `the same'.

\definition{Definition}
Two ordered diagrams with multiplicity, say $(E,r,s,f)$
and $(E',r',s',f')$, are {\it order equivalent\/} if
there is an order preserving bijection
$\Phi \: E \rightarrow E'$ so that
$$
r(e) = r'(\Phi(e)), s(e) = s'(\Phi(e)), \text{ and }
f(e) = f'(\Phi(e)).$$
We write this as $(E,r,s,f) \cong^{ord} (E',r',s',f')$.
\enddefinition

Given two order equivalent diagrams with multiplicity,
they both induce the same embedding,
providing we use the same algebras as domain for both
embeddings.

\proclaim{Theorem 6}  An embedding
 $\phi \colon \oplus_{i=1}^a T_{m_i}
        \rightarrow \oplus_{j=1}^b T_{n_j}$
is order preserving if, and only if,
 there is an ordered diagram with multiplicity
whose associated embedding is $\phi$.
\endproclaim

\demo{Proof} One direction is obvious.  For the other
we must show that $\phi$ is a direct sum of embeddings, each of which
is essentially a refinement embedding of one summand of the domain
of $\phi$.  Furthermore, each direct summand of $\phi$ must be supported
on a projection which is an interval from the lattice of invariant
projections for the co-domain of $\phi$.  Also, it is clearly sufficient
to prove the theorem for the case in which the co-domain of $\phi$
consists of a single full upper triangular matrix algebra.
\par
Thus, we assume that $\phi \colon \oplus^a_{k=1}
T_{m_k} \rightarrow T_n$
is order preserving.
For each $k$, let $\{e^k_{ij}\}$ be a matrix unit system for
$T_{m_k}$ and let $\{f_{ij}\}$ be a matrix unit system for $T_n$.
We also identify each $e^k_{ij}$ with the obvious matrix unit
in $\oplus^a_{k=1}T_{m_k}$.
For clarity, let $e^k_i$ denote the minimal diagonal projection
$e^k_{ii}$ and, similarly, let $f_i$ denote $f_{ii}$.  In each
case, the usual order on the index set corresponds to the
diagonal order on the minimal diagonal projections.
\par
The following observation is critical: we cannot have
three minimal diagonal projections, $f_b \prec f_c
\prec f_d$ and unequal integers $k$ and $l$ such that
$f_c$ is subordinate to $\phi(e^k_p)$, for some $p$;
$f_b$ is subordinate to $\phi(e^l_n)$, for some $n$;
$f_d$ is subordinate to $\phi(e^l_m)$, for some $m$;
and conjugation by $\phi(e^l_{nm})$ carries $f_d$ to $f_b$.
The reason is that if $k \ne l$ then $e^l_{nm} + e^k_p$ is
necessarily order preserving in $\oplus^a_{k=1}T_{m_k}$,
but $\phi(e^l_{nm} + e^k_p) = \phi(e^l_{nm}) + \phi(e^k_p)$
is not.  (Conjugation by the latter partial isometry
 maps $f_c$
to $f_c$ and $f_d$ to $f_b$; but $f_c \prec  f_d$ and
$f_b \prec f_c$.)
\par
Now consider $f_1$, the first diagonal projection in $T_n$.
There is a unique index $k$ such that $f_1$ is subordinate
to $\phi(e^k_1)$.  Let $1^k$ denote the projection
$e^k_1 + \dots + e^k_{m_k}$.  This operator is the projection
in the domain algebra for $\phi$ corresponding to the summand
$T_{m_k}$.  Let $\psi$ denote the mapping obtained by restricting
$\phi$ to $T_{m_k}$ and also compressing to $\phi(1^k)$.  Let $s$
be the number of minimal diagonal projections $f_i$ which are
subordinate to $\phi(1^k)$.  We retain the diagonal order on these
projections inherited from $T_n$; with respect to this order, $\psi$
is an order preserving embedding from $T_{m_k}$ to $T_s$.  By theorem~5,
$\psi$ is a direct sum of refinement embeddings.
\par
Let $\rho$ be the first summand of $\psi$ (the one for which
$\rho(e^k_1)$ contains $f_1$ as a subordinate).  Let $t$ be the
multiplicity of $\rho$.  Observe that each subordinate of
$\rho(1^k)$ precedes all of the other subordinates of $\phi(1^k)$.
This fact, combined with the critical observation above, implies
that the subordinates of $\rho(1^k)$ are $f_1, f_2, \dots, f_{tm_k}$.
In other words, $\rho(1^k)$ is an interval from the nest associated
with $T_n$ (and a leading interval, at that).
\par
It is now clear that we can split $\rho$ off from $\phi$ as a direct
summand and apply induction to what remains to see that $\phi$ must
have the desired form. \qed
\enddemo

\head 3 The Spectrum for Locally Order Preserving Embeddings
\endhead

Having described embeddings that are locally order preserving,
it is natural to consider the algebras $\indlimit(T_{n_i},\alpha_i)$ where
each $\alpha_i \: T_{n_i} \rightarrow T_{n_{i+1}}$ is locally order preserving.
There is also a smaller class of algebras,
properly contained in those with each $\alpha_i$ locally order preserving
and properly containing those with each $\alpha_i$ order preserving.
This class consists of all algebras $\indlimit(A_i,\alpha_i)$
where for each $i$ and $j$ with $i < j$, we have $\alpha_{i,j}$
is locally order preserving, where $\alpha_{i,j}$ is the composition
$$ \alpha_{j-1} \circ \alpha_{j-2} \circ
         \cdots \circ \alpha_{i+1} \circ \alpha_i.$$
We will later show (in Theorem~18)
that this class is precisely all those
strongly maximal TAF algebras
where the order preserving normalizer generates the algebra.
(The order preserving normalizer
in a TAF algebra is defined in Section~5.)

The classification in Section~8 of limit algebras with order
preserving presentations makes critical use of an invariant,
sometimes called the {\it fundamental relation}
but which we shall call the {\it spectrum}, for
subalgebras of AF \cstar-algebras which contain a canonical
masa.  We introduce the latter term because this invariant
plays a role analogous to the role played by the
spectrum (i.e.~the maximal ideal space) of an abelian
\cstar-algebra.
  This invariant was first described in the form in which
we need it in \cite{Po1}. We shall describe it briefly; a more
complete account may be found in \cite{Po4}.
\par
Let $\sA$ be a TAF algebra with diagonal $\sD$ and let $X$
be the maximal ideal space for $\sD$.  The spectrum
for $\sA$ is a topological binary relation, denoted by $R(\sA)$,
on $X$.  This relation is determined by the normalizing partial
isometries in $\sA$.  Since $\sD$ is isomorphic to $C(X)$ and since
each normalizing partial isometry acts by conjugation on $\sD$,
each normalizing partial isometry induces in a natural way
a partial homeomorphism on $X$.  (Partial, because the domain for
the homeomorphism is the subset of X which corresponds to the
initial projection of the partial isometry.)  The spectrum
 is the union of the graphs of all these partial
homeomorphisms; the topology is generated by taking
each such graph as an open and closed subset of $R(\sA)$.
\par
The spectrum can also be described in the language
of groupoids.  The enveloping \cstar-algebra for $\sA$ is a groupoid
\cstar-algebra; the spectrum for $\sA$ is the support
subsemigroupoid for the algebra $\sA$.  The main significance of the
spectrum for us is that it is a complete isometric
isomorphism invariant for triangular subalgebras of AF \cstar-algebras
when the diagonal algebras are regular canonical masas.  Effective use
of the spectrum often requires calculating a specific
representation for the spectrum.  We do this in our context
in this and the following section.

\subhead Locally Order Preserving Embeddings \endsubhead
Many of the spectra which have been described explicitly
in the literature have a common form.
Here we show this common form is precisely equivalent to the existence
of a presentation with locally order preserving embeddings.

Consider a system with locally order preserving embeddings
$$
\system {k_1} {k_1k_2} {k_1k_2k_3} {\sA}
$$
\pagebreak

Let $n_m$ denote $ k_1k_2\cdot \dots \cdot k_m $.
By Lemma~2 we obtain, for each $m \in \nbb$:
\roster
\item"{A})" A bijection between $\braks {k_1} \times
     \dots \times \braks {k_m}$ and the minimal diagonal projections of
     $D_{n_m}$ --- we will define $X_m =
     \braks{k_1} \times \dots \times \braks {k_m}$ as the
     index set for the minimal diagonal projections of $D_{n_m}$.
\item"{B})" A total order on $X_m$ (which we denote by
   $\preceq_m$) so that the bijection in A) is an order isomorphism.
   (The order on the minimal diagonal projections of $D_{n_m}$
   is the diagonal order.)
\endroster
Furthermore, the indexing and order satisfy:
\roster
\item"{a)}" $\phi_m\colon T_{n_m} \to T_{n_{m+1}}$ is given
on matrix units by the formula
$$
\phi_m(e_{(x_1, \dots ,x_m),(y_1,\ldots,y_m)}) =
 \sum_{j=1}^{k_{m+1}} e_{(x_1, \dots ,x_m,j),(y_1,\ldots,y_m,j)}.
$$
\item"{b)}" If $i<j$ then $(x_1, \dots ,x_m, i) \preceq_{m+1}
 (x_1, \dots ,x_m, j)$.
\item"{c)}" If $(x_1, \dots, x_m) \preceq_m (y_1, \dots,y_m)$
then $(x_1, \dots ,x_m, j) \preceq_{m+1}  (y_1, \dots ,y_m, j)$.
\endroster

\definition{Definition}
A sequence of orders $\preceq_m$
on the sets $X_m$ satisfying properties b) and c) is
said to be {\it coherent\/}.
\enddefinition

Let $\displaystyle X = \prod_{j=1}^{\infty} \braks{k_j}$
and give $X$ the product topology.
Then $X$ is isomorphic to the maximal ideal space of $\sD$, the diagonal
of $\sA$, where $\sD = \varinjlim D_{n_m}$.
Let $R(\sA)$ denote the spectrum, a topological binary
relation on $X$.
\par
It is routine to show that the indexing above yields:
$$ \split xR(\sA)y \Longleftrightarrow
  & \text{ There exists } m \in \nbb \text{ such that } x_n = y_n
       \text{ for all } n \ge m \\
  & \text{ and } (x_1, \dots ,x_m) \preceq_m (y_1, \dots ,y_m).
\endsplit $$
While there are many choices for $m$,
coherence guarantees that the initial segments are ordered the same way
for any choice of $m$ giving common tails.
\par
More colloquially,
$$ \split xR(\sA)y \Longleftrightarrow
  & \text{$x$ and $y$ have the same tails and the initial segments} \\
  & \text{are ordered with respect to a coherent sequence of orders.}
\endsplit $$

Conversely, let $\displaystyle X = \prod_{j=1}^{\infty} \braks{k_j}$, \
$\displaystyle X_m = \prod_{j=1}^m \braks {k_j}$, \     and
$\preceq_m$ be a total order on $X_m$.
Assume that the sequence of orders is coherent.
\par
Let $R$ be a topological binary relation on $X$
defined by:
$$ \split xRy \Longleftrightarrow
  & \text{ There exists } m \in \nbb \text{ such that } x_n = y_n
        \text{ for all } n > m \\
  & \text{ and } (x_1, \dots , x_m) \preceq_m (y_1, \dots ,y_m).
\endsplit
$$
The topology is given by taking  (for each $m \in \nbb$
and $a,b \in X_m$) the following sets as a basis
of clopen sets:
$$
E_{a,b} = \{\,(x,y)\in X \times X \colon
  x_n = y_n \text{ for } n>m, \,x_n = a_n \text{ and }
y_n = b_n \text{ for } 1 \le n \le m\,\}
$$
It is clear that each $E_{a,b} \subset R$.

\definition{Definition}
A topological binary relation is
{\it coherent\/} if it is isomorphic (as a topological
binary relation) to one of the form above.
Actually, the form given above is only one representation
of the topological binary relation, so it would be more precise
to say that it is a topological binary relation {\it with a coherent
representation}.
For the sake of brevity, we use the shorter term.
\enddefinition

\proclaim{Theorem 7}
If $\sA$ is the direct limit of a system,
$$
\system {k_1} {k_1k_2} {k_1k_2k_3} \sA
$$
where $\phi_i$ is locally order preserving for each $i$,
then $R(\sA)$ is coherent.

Conversely, if $R$ is coherent, then $R \cong R(\sA)$
where $\sA$ is a direct limit of  such a system.
\endproclaim

\demo{Proof}
We have already proved the first statement in the theorem,
so only the converse remains to be proved.

As usual, $n_m = k_1 \cdot \dots \cdot k_m$
and $X_m =\braks {k_1} \times \dots \times \braks {k_m}$.
Assume $R$ is isomorphic to a relation with the
form above.   Use $X_m$
as the index set for the minimal diagonal
 projections in $D_{n_m}$; order the
minimal diagonal projections according to the order on $X_m$; let $T_{n_m}$
be the algebra of upper triangular matrices with respect to this order.
Define $\phi_m\colon D_{n_m} \to D_{n_{m+1}}$ by
$$
\phi_m(e_{(x_1, \dots ,x_m)}) =
 \sum_{j=1}^{k_{m+1}} e_{(x_1, \dots ,x_m, j)}
$$
and extend $\phi_m$ to a locally order preserving embedding
$T_{n_m} \to T_{n_{m+1}}$.  The coherence guarantees
that $\phi_m(T_{n_m}) \subseteq T_{n_{m+1}}$.
\par
Let $\sA = \varinjlim (T_{n_m}, \phi_m)$.
Then it is clear that $R \cong R(\sA)$.
\qed
\enddemo

\subhead An Intermediate Family \endsubhead
There is an analogue of Theorem~7 for the class of algebras
$\indlimit(A_i,\alpha_i)$ where each $A_i$
 is the upper triangular matrices
of some full matrix algebra and each composition of embeddings
$\alpha_{i,j}=\alpha_{j-1} \circ \cdots \circ \alpha_i$ is locally order
 preserving.

We may rephrase the second line of (3) in Section 2  as follows:
the diagonal order restricted to $\{ (i,l) \stbar l \in [k] \}$ induces
the same order on $[k]$, for each choice of $i \in [n]$.
The appropriate generalization is that given $i$ and $j$ with $i < j$, then
the diagonal order on $X_j$ restricted to
$$\{ (x_1,\ldots,x_i,a_{i+1},\ldots,a_j) \stbar a_l \in [k_l]
\text{ for all $l$ with } i < l \le j \}$$
induces the same order on
$[k_{i+1}] \times [k_{i+2}] \times \cdots \times [k_j]$
for each choice of an element $(x_1,\ldots,x_i)$ in $X_i$.

\definition{Definition}
A spectrum with a coherent representation that satisfies
this additional condition for all $i$ and $j$ with $i < j$
will be called {\it hypercoherent}.
\enddefinition

\proclaim{Theorem 8}
If $\sA$ is the direct limit of a system,
$$
\system {k_1} {k_1k_2} {k_1k_2k_3} \sA
$$
where for all $i$ and $j$ with $i < j$, the composition
$\phi_{i,j} = \phi_{j-1} \circ \cdots \circ \phi_i$ is locally order
preserving, then $R(\sA)$ is hypercoherent.

Conversely, if $R$ is hypercoherent, then $R \cong R(\sA)$
where $\sA$ is a direct limit of such a system.
\endproclaim

\demo{Proof}
The following observation is the key to proving both directions:
given a matrix unit $e = e_{(x_1,\ldots,x_i),(y_1,\ldots,y_i)}
 \in T_{k_1 \cdots  k_i}$,
then $ \phi_{i,j}(e)$ is order preserving in $T_{k_1 \cdots k_j}$
if and only if the diagonal order on $X_j$ restricted to
$$\{ (x_1,\ldots,x_i,a_{i+1},\ldots,a_j) \stbar a_l \in [k_l]
\text{ for all $l$ with } i < l \le j \}$$
and
$$\{ (y_1,\ldots,y_i,a_{i+1},\ldots,a_j) \stbar a_l \in [k_l]
\text{ for all $l$ with } i < l \le j \}$$
both induce the same order on $[k_{i+1}] \times [k_{i+2}] \times \cdots \times
 [k_j]$.

If $\sA$ is the direct limit of a system with each $\phi_{i,j}$ locally
order preserving, then
the observation immediately implies that $R(\sA)$ is hypercoherent.
Conversely, given a hypercoherent spectrum $R$,  by
Theorem~7 we can construct a system
$$
\system {k_1} {k_1k_2} {k_1k_2k_3} \sA
$$
with locally order preserving embeddings such that $R \cong R(\sA)$.
As $R(\sA)$ is hypercoherent, the observation
 implies that for each $i$ and $j$
with $i < j$, and for each matrix unit $e \in T_{k_1\cdots k_i}$,
the image $\phi_{i,j}(e)$ is order preserving in $T_{k_1\cdots k_j}$.
Thus, each $\phi_{i,j}$ is locally order preserving, as required.
\qed
\enddemo

\head 4 Spectra for Order Preserving Embeddings
\endhead

We turn now to the spectrum for algebras
$\indlimit(T_{n_i}, \alpha_i)$ where each
$\alpha_i \: T_{n_i} \rightarrow T_{n_{i+1}}$ is order preserving.
After describing these spectra in the first subsection,
we characterize the gap points in Theorem~9 and
then, in Theorem~13, prove the
first step in our classification of such algebras.
The last two theorems of the section show that these algebras
are analytic as are the algebras $\indlimit(A_i,\alpha_i)$,
where the $A_i$ are allowed to be direct sums of $T_n$'s
and the $\alpha_i$ are still order preserving.

\def\r#1#2{r^{(#1)}_{#2}}
\def\F#1#2{F^{(#1)}_{#2}}
\def\orb#1{\Cal O(#1)}
\def\clorb#1{\overline{\orb{#1}}}
\def\row#1#2{(#1_1, \dots , #1_{#2})}

\def\sn#1{\text{sn}(#1)}
%
%
\subhead The Spectrum \endsubhead
Let, for each $n$,
$\phi_{\r n{}} = \rho_{\r n1} \oplus \dots \oplus \rho_{\r n{t_n}}$
where $\r n{}$ is the $t_n$-tuple $(\r n1, \dots , \r n{t_n})$
and consider the direct system:
\setbox1 = \hbox{$\scriptstyle \phi_r(1)$}
$$T_1 @>\phi_{\r 1{}}>> T_{k_1} @>\phi_{\r 2{}}>>
      T_{k_1k_2} @>\phi_{\r 3{}}>> \dots
 @>\hphantom{\phi_{\r 3{}}}>> \sA. $$
Since direct sums of refinement embeddings are order preserving,
we can apply the results of the previous section.
In this case, we describe explicitly the hypercoherent orderings
which extend naturally
the lexicographic orders used for refinement embeddings and
the reverse lexicographic orders used for standard embeddings.
This description also subsumes that given in \cite{HPo} for the
spectrum of an alternation algebra.

To fix notation, note that the
 maximal ideal space of $\sA$ is isomorphic to
$\displaystyle X = \prod_{n=1}^{\infty} \braks {k_n}$
where $k_n = \r n1 + \dots + \r n{t_n}$.
Let $\F n1, \dots , \F n{t_n}$ be the partition of $\braks {k_n}$
corresponding to $\r n1, \dots ,\r n{t_n}$.
Specifically,
\newcount\fsetpage
\fsetpage = \the\count0
$$\align
 \F n1 & = \{\,1, \dots, \r n1\,\},\\
 \F n2 & = \{\, \r n1 +1, \dots ,\r n1 + \r n2\,\},\\
 \vdots& \\
 \F n{t_n} & = \{\, \r n1 + \dots + \r n{t_n-1} + 1 , \dots , k_n \,\}.
\endalign $$
For each integer $x \in \braks {k_n}$, define
$i_n(x)$ to be the unique index~$s$ so that $x \in \F ns$.

Applying Lemma~2 to
 $\phi_{r^{(1)}}\colon T_{k_1} \to T_{k_1k_2}$
gives the following ordering on
 $\braks {k_1} \times \braks {k_2}$:
$$
(x_1, x_2) \preceq_2 (y_1, y_2) \qquad\text{if }
\left\{
\aligned
& i_2(x_2) < i_2(y_2), \text{ or} \\
& i_2(x_2) = i_2(y_2) \text{ and }   x_1  <   y_1 , \text{ or} \\
& i_2(x_2) = i_2(y_2), \ x_1 = y_1, \text{ and }
       x_2 \le y_2. \\
\endaligned \right.
$$
\par
Next, define a total order on each of the sets
$\braks{k_1} \times \dots \times \braks {k_n}$ recursively:
\roster
\item The order $\preceq_2$ is defined on $\braks {k_1}
      \times \braks {k_2}$ as above.
\item The order $\preceq_n$ is defined on
      $\braks{k_1} \times \dots \times \braks {k_n}$ by
      treating it as the Cartesian product of
      $\braks {k_1} \times \dots \times \braks {k_{n-1}}$ carrying
      the order $\preceq_{n-1}$ and $\braks {k_n}$ with the usual total
      order and again applying the procedure above.
\endroster

By Theorem 7, $xR(\sA)y$ if, and only if, there exists
$ m \in \nbb $ such that $ x_n = y_n $ for all $ n > m $
and $ (x_1, \dots , x_m) \preceq_m (y_1, \dots ,y_m)$.
In terms of the $i_n$ functions, if $xR(\sA)y$, then
\roster
\item $x=y$, or
\item There is a $q$ such that $i_q(x_q) <
   i_q(y_q)$ and $i_n(x_n) = i_n(y_n)$ for
   all $n>q$, or,
\item $i_n(x_n)=i_n(y_n)$ for all $n$, and
  there is a $q$ such that $x_q < y_q$ and $x_m = y_m$
  for all $m<q$.
\endroster
Informally, to determine whether $xR(\sA)y$ or $yR(\sA)x$
for $x$ and $y$ with the same tails we compare
the initial segments of $x$ and $y$.
First we look for the highest coordinates which
belong in different $\F n{}$-sets;
if this does not occur, we look for the lowest coordinates which differ.
The order of the $\F n{}$-sets or of the coordinates themselves,
as appropriate, determines the order between $x$ and $y$.
In short, the order is reverse lexicographical for $\F n{}$-sets,
then lexicographical for the coordinates themselves.

\subhead Gap Points \endsubhead
The material in this section on gap points and the next on first
refinement multiplicities follows the line of argument in \cite{HPo},
there given for the special case of alternation algebras.

The description of the gap points is easy; the verification
that the description is correct is tedious.
\par
Let $\orb x = \{\,z \stbar  zR(\sA)x\,\}$ be the orbit of $x$
and $\clorb x$ the closure of the orbit.
\definition{Definition}
 We define $x$ to be a {\it gap point\/}
 if there is a point $y$  such that
 $y \notin \clorb x$ and $\clorb y = \clorb x \cup \{\,y\,\}$.
\enddefinition

\remark{Remark}
It would be  more accurate to call
$x$ a {\it left gap point\/}.  Then $y$ is the corresponding
{\it right gap point\/} and could well be denoted by $x^+$.
\endremark

There is one possible exception to the characterization
of gap points in the following theorem.  Let $x^{\infty}$
denote the sequence $(k_n)$.  Thus, each coordinate of
$x^{\infty}$ is the maximal element of the $\F n{}$-set
with maximal index.  This point is exceptional for the
condition below only in the case of a refinement algebra.

\proclaim{Theorem 9}
A point $x \ne x^{\infty}$ is a gap point if,
and only if, there is an integer $p$ so that for all
$n>p$, $i_n(x_n) = \max \F n1$.  In other words, for
large $n$, $x_n$ is the largest element of the first
$\F n{}$-set; viz., $x_n = \r n1$.
\endproclaim
As we shall see, the gap points effectively determine
up to a finite factor
the supernatural number of the sequence $(\r n1)$ of
multiplicities of the first refinement summands.
\par
We shall need several preliminary facts about orbits
in order to prove the theorem.
For a point $y \in X$, define
$$
W_p(y) = \{\,z \in X \mid z_1 = y_1, \dots , z_p = y_p\,\}.
$$
When $y$ is clearly understood, we write $W_p$ instead.
 These sets are
open, as well as closed, and form a basis for the topology
at $y$.  Convergence in the topology is pointwise:
$$
z(n) \to z \Longleftrightarrow z(n)_i \to z_i
\text{ as } n \to \infty, \text{ for each } i.
$$
Lemma~10 below, together with
the observation that if $x$ is a gap point,
then $\orb x$ is not dense in $X$,
will establish the necessity of the condition
$x_n \in \F n1$ for all large $n$.
\proclaim{Lemma 10}
Assume that there are infinitely many
$n$ with $x_n \in \F np$ for $p>1$.  Then $\orb x$
is dense in $X$.
\endproclaim

\demo{Proof}
Let $y \in X$ be arbitrary.
We need merely show that, for each $p>0$,
$W_p$ contains a point in the orbit $\orb x$.
Given $p$, choose $n>p$ such that $x_n \notin \F n1$.
Let $z_n$ be any element in $\F n1$.
For $1 \le t \le p$, let $z_t = y_t$.
For all other $t\ne n$, let $z_t = x_t$.
Then $z = (z_1, z_2, \dotsc )$ lies in $ W_p \cap \orb x$.
Thus, $\orb x$ is dense in $X$.
\qed
\enddemo

Let
$$ \aligned
A_p &= \{\,x \in X \mid x_t \le \r t1 \text{ for }
    t \ge p\,\} \\
&=\{\,x \in X \mid i_t(x_t) =1 \text{ for } t \ge p\,\} \\
&=\prod_{n=1}^{p-1}   \braks {k_n} \times
  \prod_{n=p}^{\infty} \F n1
\endaligned $$
Observe that each $A_p$ is a closed set.
\par
Now suppose that $y \in X$ and $\orb y$ is not dense in $X$.
Then there exists an integer $p$ such that
 $y_n \in \F n1$, for all $n \ge p$.
If $x \in \orb y$, then $i_n(x_n) =1$ for all $n \ge p$
also, so $\orb y \subseteq A_p$.  Since $A_p$ is closed,
$\clorb y  \subseteq A_p$.

With $y$ as above, suppose that $x \in A_p$ and that
$\row xq \preceq \row yq$ for some $q \ge p$.  Then
$\row xn \preceq \row yn$ for all $n$ satisfying
$p \le n \le q$.  This follows immediately from  the
definition of $\preceq$ and the fact that
$i_n(x_n) = 1 = i_n(y_n)$ for $p \le n \le q$.
If we actually have $\row xq \prec \row yq$, then
we also have $\row xn \prec \row yn$ for all $n \ge q$.

With $y$  still as above and $x$ in $\orb y$, there
is an integer $q \ge p$ such that $y_n = x_n$ for all
$n \ge q$.  If $y \ne x$, then
$\row xq \prec \row yq$.  Consequently,
$$
x \in \orb y \Longrightarrow \row xn
\preceq \row yn \text{ for all } n \ge p.
$$
\proclaim{Lemma 11}
Assume that $\orb y$ is not dense
in $X$, i.e. that there is an integer $p$ such that
$i_n(y_n) =1$, for all $n \ge p$.  Then
$$
x \in \clorb y \Longleftrightarrow
  \text{there exists } s \ge p \text{ such that }
 \row xn \preceq \row yn, \text{ for all } n \ge s.
$$
\endproclaim

\demo{Proof}
Let $x \in \clorb y$.  Let $n \ge p$ and let $W_n$ be
the corresponding neighborhood of $x$.  There is
an element $z \in \orb y$  such that $z \in W_n$.
Consequently, $\row xn = \row zn \preceq \row yn$.
This establishes the implication $\Rightarrow$.

For the converse, we may assume, based on the
remarks above, that $\row xn \preceq \row yn$
for all $n \ge p$.  Fix $n \ge p$ and define a
point $z \in X$ by
$$
z_t = \left\{
\aligned
&y_t \qquad \text{if } 1 \le t \le n,  \\
&x_t \qquad \text{if } n < t.
\endaligned  \right.
$$
Then $z \in \orb y$ and $z \in W_n$.  Thus, every
neighborhood of $x$ intersects $\orb y$, hence
$x \in \clorb y$.
\qed
\enddemo

\proclaim{Corollary 12}
  If $x \in \clorb y$ then $\clorb x \subseteq \clorb y$.
\endproclaim

\demo{Proof} Apply Lemma~11. \qed \enddemo

We are now ready to prove Theorem~9.
\demo{Proof of Theorem~9}
Let $x$ be a (left) gap point.  Lemma~10 shows that there
is an integer $p$ such that $i_n(x_n) =1$ for all $n \ge p$.
We need to show further that $x_n = \max \F n1$ for all but
finitely many $n$.
\par
Assume the contrary; that is, assume that $i_n(x_n) = 1$
for all $n \ge p$ but that $x_n \ne \max \F n1$
for infinitely many values of $n$.  We must show that
there is no point $y$ which satisfies $y \notin \clorb x$
and $\clorb y = \clorb x \cup \{y\}$.  Clearly, $y$ cannot
satisfy these two properties if $\orb y$ is dense in $X$ or
if $\clorb x \not \subseteq \clorb y$.  So we may reduce
to the case in which    $y$ satisfies:
             $x \in \clorb y$,
             $x \ne y$
and           there is an integer $q$ such that
 $i_n(y_n) =1$ for all $n \ge q$.
Without loss of generality, we may assume that $q \ge p$.
\par
By Lemma~11 there is an integer $s \ge q$ such that
$$
\row xn \prec \row yn \text{ for all } n \ge s.
$$
        Let $m$ be any integer such that $m > s$ and
$x_m \ne \max \F m1$.  Define $z = z(m)$ by
$$
z_t = \left\{
\aligned
&x_t \\ &x_m +1
\endaligned \right. \qquad
\aligned
&\text{if } t \ne m, \\
&\text{if } t = m.     \\
\endaligned
$$
\par
It is evident from Lemma~11 that $z \notin \clorb x$.
Lemma~11 also shows that $z \in \clorb y$.  Indeed,
$i_t(z_t) = i_t(x_t) = 1 = i_t(y_t)$ for all $t > s$ and
$\row zs \prec \row ys$.  (Note: this uses $i_m(z_m) =
i_m(x_m+1) = 1$.)  Consequently, $\row zn \prec \row yn$
for all $n\ge s$ and $z \in \clorb y$, as desired.
\par
Since there are infinitely many $m$ such that
$x_m \ne \max \F m1$, there are infinitely many
distinct points $z(m)$ which are in $\clorb y$
but not in $\clorb x$.  Thus, in this case,
$\clorb y \ne \clorb x \cup \{y\}$.
\par
Next, we prove the converse.  So assume that $x$
satisfies $x_n = \max \F n1$ for all $n \ge p$.
We must produce an element $y \in X$ such that
$y \notin \clorb x$ and $\clorb y = \clorb x \cup \{y\}$.
We consider separately two cases.

{\it Case 1.\/} There is an integer $m$ such that
$x_m \ne \max \F mj$ and $x_n = \max \F nj$ for $n > m$.
(Here, $j = i_m(x_m)$ or $j = i_n(x_n)$ as appropriate.  Of
course, if $n \ge p$  then $j=1$.   Also, note that
$x_t \ne \max \F tj$ is possible for only finitely
many $t$, so there must be a maximal such $t$ if there are
any.)
\par
Define $y \in X$ by:
$$
y_t = \left\{
\aligned
&x_t,\\  &x_m +1, \\ &\min \F tj, \\
\endaligned \right. \qquad
\aligned
&\text{if } 1 \le t \le m-1,  \\
&\text{if } t=m, \\
&\text{if } t > m, \quad j = i_t(x_t).
\endaligned
$$
We have $i_n(y_n) = i_n(x_n)$ for all $n$, and
in particular, $i_n(y_n) = 1$ for all $n \ge p$.
 From this and the fact that $x_m < y_m$, it is clear
that $\row xn \prec \row yn$ for all $n \ge m$.
So $x \in \clorb y$ and hence (by Corollary~12)
$\clorb x \subseteq \clorb y$.
It is also clear that $y \notin \clorb x$.
\par
It remains to show that if $z \in \clorb y$ and $z \ne y$,
then $z \in \clorb x$.  Given such a $z$, let $r >m$ be such
that $\row zn \prec \row yn$ for all $n \ge r$.  It suffices to
show that $\row zn \preceq \row xn$ for all $n \ge r$.
Fix $n \ge r$.
If $i_t(z_t) < i_t(y_t)$ for some $t \le n$ (a possibility
only if $t <p$), then $\row zn \prec \row xn$.
If $i_t(z_t) = i_t(y_t)$ for all $t$, then there is an
integer $q$ such that $z_q < y_q$ and $z_t = y_t$ for
all $t < q$.  We cannot have $q > m$, since
$y_t = \min \F tj$ for $t > m$.  If $q<m$, the
$\row zn \prec \row xn$ as desired.  If $q = m$
then $z_m \le x_m$.  If $z_m < x_m$, we again
have $\row zn \prec \row xn$.
Finally, if $z_m = x_m$  then, since
$i_t(x_t) = \max \F tj$ for all $t>m$ and
$i_t(z_t) = i_t(y_t) = i_t(x_t)$ for all $t$,
we have $\row zn \preceq \row xn$.
This exhausts all possibilities and case 1 is complete.

\resultspace
{\it Case 2.\/}  For all $n$, $x_n = \max \F nj$, where
 $j = i_n(x_n)$.  Let $q$ be the least integer such
that $x_q \ne k_q$, i.e., $i_q(x_q)$ is not maximal
in the set of indices for $\F q{}$.
There must be such an integer $q$, since we assume
that $x$ is not the exceptional point $x^{\infty}$.
(This could, in fact, happen only if the algebra is
actually a refinement algebra: $\F n1 = \braks {k_n}$
for all large $n$.)  We now define $y$ by:
$$
y_t = \left \{
\aligned
&1, \\  &x_q +1, \\ &\min \F tj, \\
\endaligned \right.  \qquad
\aligned
&\text{if } 1 \le t \le q-1, \\
&\text{if } t = q, \\
&\text{if } t>q \text{ and } j= i_t(x_t). \\
\endaligned
$$
Observe that $i_t(y_t) = 1$ if $1 \le t \le q-1$,
that $i_q(y_q) = i_q(x_q) + 1$, and that
$i_t(y_t) = i_t(x_t)$ for all $t > q$.
In particular, $y_t = \min \F tj$ for all $t$.
\par
For $n \ge q$ it is clear that $\row xn \prec \row yn$.
This implies that $y \notin \clorb x$ and $x \in \clorb y$,
which, by Corollary~12, implies that $\clorb x \subset \clorb y$.
It remains to prove that $\clorb y = \clorb x \cup \{y\}$.
\par
Let $z \in \clorb y$ be such that $z \ne y$.  We must have
$i_t(z_t) = 1 = i_t(y_t)$ for all large~$t$.  We claim that
$i_t(z_t) \ne i_t(y_t)$ for some~$t$.  Indeed, if
$i_t (z_t) = i_t(y_t)$ for all~$t$, then the two facts
$\row zn \preceq \row yn$ for all large $n$ and
$y_t = \min \F tj$ for all $t$ imply that
$\row zn = \row yn$ for all large $n$.  But this means that
$z=y$.
\par
Thus, there is an integer $s$ such that $i_s(z_s) < i_s(y_s)$
and $i_t(z_t) = i_t(y_t)$ for all $t > s$.  Note that
$s <\max\{\,p,q\,\}$, since if $t \ge \max \{\,p,q\,\}$,
then $i_t(y_t) = i_t(x_t) = 1$.
If $s >q$ then $i_s(z_s) < i_s(y_s) = i_s(x_s)$
and $i_t(z_t) = i_t(y_t) = i_t(x_t)$ for all $t >  s$.
Hence, $\row zn \prec \row xn$ for all $n \ge s$ and
$z \in \clorb x$.
We cannot have $s < q$, since if
we do,  $i_s(y_s) = 1$, which is incompatible with
$i_s(z_s) < i_s(y_s)$.  So we are left with the case in
which $s = q$.  We then have $i_s(z_s) \le i_s(x_s)$.
For $t < s$, $i_t(x_t)$ is maximal, so it follows that
$i_t(z_t) \le i_t(x_t)$ for all $t$.
Since $x_t = \max \F tj$ for all $t$, we have
$\row zn \preceq \row xn$ for all $n$.

Thus in all cases, $z \in \clorb x$.
This completes the proof of the theorem.
\qed
\enddemo

\subhead First Summand Refinement Multiplicities \endsubhead
Let $x$ be a gap point and let $y \in \orb x$.  Let
$p$ be an integer such that $i_n(x_n) = 1$ for all
$n \ge p$ and $x_n = \max \F n1$ for all $n \ge p$.
We may assume that $p$ is the least integer with
these properties, though it is not actually necessary
to do so. Lemma~11 implies
$$
\clorb x = \{\,z \colon  \row z{p-1} \preceq  \row x{p-1}
   \text{ and } z_n \in \F n1 \text{ for } n \ge p\,\}.
$$
If we let $h$ denote the number of elements of
$\braks {k_1} \times \cdots \times \braks {k_{p-1}}$ which precede
$\row x{p-1}$, then we have
$$
\clorb x \cong \braks h \times \prod_{n = p}^{\infty}
   \braks {\r n1}.
$$
We endow $\clorb x$ with the relative topology
and order that it inherits from $R(\sA)$ and
$ \braks h \times \prod_{n = p}^{\infty}
   \braks {\r n1}$
with the usual spectrum associated with
a refinement algebra with supernatural number
$h\prod_{n=p}^{\infty}\r n1$.

\proclaim{Theorem 13}
Suppose  that $\sA$ and $\sB$ are two direct limits of
systems where the embeddings are
order preserving and hence are
 direct sums of refinement
embeddings.  Assume that $\sA$ and $\sB$ are
isometrically isomorphic.
Let $(\r n1)$ be the sequence of multiplicities of the
first refinement summand in the embeddings for $\sA$.
Let $(s_1^{(n)})$ be the corresponding sequence for $\sB$.
Let $\sn r$ and $\sn s$ be the supernatural numbers for
these two sequences.  Then there are finite numbers
$a$ and $b$ such that $a \, \sn r = b \, \sn s$.
\endproclaim

\demo{Proof}
Let $x$ be the point with $x_n = \max \F n1$ for all  $n$.
We assume that $\sA$ and $\sB$ are not refinement algebras,
since a stronger result is known in that case ($\sn r = \sn s$).
So $x$ is not an exceptional point.  Let $\beta$ be a
spectrum isomorphism of $R(\sA)$ onto $R(\sB)$.
Then $\beta (x)$ is a gap point in $R(\sB)$ and
$\beta$ restricted to $\clorb x$ is a spectrum
isomorphism of $\clorb x$ onto $\clorb {\beta(x)}$.
The description of closed orbits above yields the theorem.
\qed
\enddemo
\par
\subhead Analyticity \endsubhead
The detailed description of the spectrum for a limit algebra
with order preserving embeddings makes it fairly simple to
prove that these algebras are all analytic.
\par
Analytic algebras have been studied in  detail in
papers such as \cite{V1,PePW2,SV,PW};
 we refer the reader to these
sources for a complete description of the notion.  A practical
working definition of analyticity can be given in terms of
the existence of a cocycle on the spectrum of the algebra.
If~$R$ is a topological binary relation, a {\it cocycle}~$c$
on $R$ is a continuous function $c\:R \longrightarrow \Bbb R$
satisfying the ``cocycle'' property: $c(x,y) + c(y,z) = c(x,z)$
for all~$x$, $y$ and $z$ such that $(x,y) \in R$ and $(y,z) \in R$.
(Technically, $c$ is a 1-cocycle.)
  In our situation, where $R$
is the spectrum of an AF algebra, $R$ will be a relation on a
compact Hausdorff topological space $X$ (the maximal ideal space
of a canonical masa).  We sometimes
 find it convenient, albeit a little
imprecise, to refer to $c$ as a cocycle on $X$.
\par
A cocycle is {\it locally constant}  if it is constant on
   some neighborhood of each point of its domain.  Locally
constant cocycles have been studied in \cite{VW,V2}.
Clearly, a locally constant function is always
 continuous.
\par
\proclaim{Theorem 14}
If $\sA = \indlimit(T_{n_i}, \alpha_i)$ is a direct limit
with every $\alpha_i$ order preserving, then $\sA$ is an
analytic algebra.  Furthermore, there is a locally constant
cocycle defined on the spectrum of $\sA$.
\endproclaim
\demo{Proof}
 We may assume that the spectrum, $R(\sA)$, of $\sA$
is a topological binary relation defined on the set
$X = \prod_{i=1}^{\infty} \braks{k_i}$
in the fashion described 
 at the beginning of this section.
Let $X_m = \prod_{i=1}^m\braks{k_i}$, for each positive
integer $m$.  The sequence of sets $X_m$ carries a coherent
family of orders $\preceq_m$ which determines the order on $X$.
In order to define a cocycle $c$ on $X$, it will suffice to
define a sequence of cocycles $c_m$ on $X_m$ which satisfy
the properties:
\roster
\item{}$c_m(x,y) \ge 0$ if, and only if $x \preceq_m y$,
for all $x,y \in X_m$.
\item{}If $x = (x_1,\dots,x_m) \in X_m$ and $j \in \braks{k_{m+1}}$,
let $(x,j)$ denote $(x_1,\dots,x_m,j)$, an element of $X_{m+1}$.
Then $c_{m+1}((x,j),(y,j)) = c_m(x,y)$,
     for all $x,y \in X_m$ and $j \in \braks{k_{m+1}}$.
\endroster
\par
Indeed, given such a sequence of cocycles,   we can define
a cocycle $c$ on $X$ as follows:
if $xR(\sA)y$ then there is an integer $p$ such that $x_i = y_i$
for all $i \ge p$.  Define $c(x,y) = c_p((x_1,\dots,x_p),(y_1,\dots,
y_p))$.  Property (2) guarantees that $c$ is well-defined.
  The construction
makes it clear that $c$ is locally constant, and hence continuous.
Property (1) ensures that $\sA$ is the analytic algebra determined by $c$.
\par
Since $X_1 = \braks {k_1}$ and $\preceq_1$ is the usual order
on integers, $c_1$ is uniquely determined by specifying the $k_1-1$
values $c(i,i+1)$, for $1 \le i \le k_1-1$.  The only constraint
imposed by the properties above is that these numbers all be
positive.  After $c_1$ has been selected, we can define the
remaining $c_m$'s recursively.  The recursive step is identical
at each stage; furthermore, a change in notation will make this
step much easier to write down.
\par
Since $X_{m-1}$ is totally ordered by $\preceq_{m-1}$, this
set is order isomorphic to $\braks {k}$ with the usual order,
where $k$ is the cardinality of $X_{m-1}$.  So in place of
$c_{m-1}$, we may assume that we have a cocycle $c_1$ defined
on $\braks{k}$ with the property that $c(i,j) \ge 0$
whenever $i \le j$.  For simplicity of notation, let $n$
denote $k_m$ and $\preceq$ denote $\preceq_m$.
 Let $Y = \braks{k}\times\braks{n}$.  If
$r_1, \dots, r_p$ are the refinement multiplicities of
the embedding associated with this step, then
$n = r_1 + \dots + r_p$ and the order on $Y$ can be described
as follows: the first elements of $Y$ are
the elements of $\braks{k}\times F_1 = \braks{k}\times\braks{r_1}$
in the lexicographic order.  Next come all the elements
of $\braks{k}\times F_2$, again in the lexicographic
order.  Continue with the groups $\braks{k}\times F_3, \dots,
\braks{k}\times F_p$, always with the lexicographic order within
each group.  The sets $F_i$ are the sets defined
on page \the\fsetpage.
\par
Our task then is to define a cocycle $c_2$ on $Y$ subject to the
properties:
\roster
\item{}$c_2(x,y) \ge 0$ if, and only if $x \preceq y$.
\item{}$c_2((i,t),(j,t)) = c_1(i,j)$ for all $i,j \in \braks{k}$
  and all $t \in \braks{n}$.
\endroster
We will define $c_2$ at all pairs $(\alpha, \beta)$ where
$\beta$ is the immediate successor of $\alpha$ in $Y$ under
the order $\preceq$.  The cocycle property then determines $c_2$.
If $\alpha$ is the last element in the group
$\braks{k}\times F_t$ and
$\beta$ is the first element in the next group, then
$c_2(\alpha,\beta)$ may be chosen arbitrarily, so long as it
is positive.  (This is the case in which $\alpha
= (k,r_1+\dots+r_t)$ and
 $\beta = (1,r_1+\dots+r_t+1)$.)
\par
Now suppose that $\alpha \in \braks{k}\times F_t$ and that
$\alpha$ is not the last element of this subset of $Y$.  There
are two cases to be distinguished.  In the first,
$\alpha = (i,j)$, where $i \in \braks{k}$ and where $j$
satisfies $r_1 + \dots + r_{t-1} +1 \le j <
r_1 + \dots + r_t$.  In this
case, the immediate successor to $\alpha$ is
$\beta = (i,j+1)$.  If we let $d = \min\{c_1(i,i+1)
\stbar 1 \le i \le k-1\}$, then we may define
$\displaystyle c_2(\alpha,\beta) = \frac d{r_t}$.
In the second case, $j = r_1 + \dots + r_t$ and
$i < k$; the immediate successor of $\alpha$ is
$\beta = (i+1, r_1 + \dots r_{t-1} +1)$.
 In this situation, we define
$\displaystyle c_2(\alpha,\beta) =  c_1(i,i+1) -
\frac {r_t - 1}{r_t}d$.
\par
We have now defined $c_2$ at all pairs $(\alpha,\beta)$,
where $\beta$ is the immediate successor of $\alpha$.
There is a unique extension of $c_2$ to a cocycle defined
on all of $Y$ and this extension satisfies the
two required properties. \qed
\enddemo

\remark{Remark}
The proof in Theorem~14 can easily be extended to a slightly
more general situation: the case in which the ``building block''
algebras are maximal triangular subalgebras of finite dimensional
\cstar-algebras; i\.e\., direct sums of $T_n$'s.
The proof of the theorem proceeds by defining the locally
constant cocycle on subsets of the spectrum which are graphs
of matrix units.  This is thinly disguised by the reductions
which were made in order to achieve notational simplification.
The specific representation of the spectrum is convenient for
expressing the proof but is not essential; use of the sets
$\braks{k_1} \times \dots \times \braks{k_n}$ as index sets
for the matrix units in corresponding finite dimensional
algebras would enable one to define the cocycle on the graphs
of the matrix units in the spectrum.
\par
  In the argument presented in the
theorem, the construction of
the cocycle~$c_2$ on the image of each refinement
embedding (more precisely, on the graphs of the matrix
units in each image) does not depend on using the same
cocycle~$c_1$ on the domain for each refinement embedding.
Consequently, if we had an order preserving embedding from
a direct sum of $T_n$'s into a single $T_n$ with a different
cocycle associated with each summand of the domain, we could
still use the same procedure to construct a cocycle for the
range.  If the range is also a direct sum of
$T_n$'s, we simply construct a cocycle for each summand of
the range and put these together.  Thus, the following is true:
\par
\proclaim{Theorem 15}
If $\sA$ is a triangular subalgebra of an AF \cstar-algebra
and if $\sA$ has a presentation with order preserving embeddings
between direct sums of $T_n$'s,
then $\sA$ is an analytic algebra with a locally constant
cocycle defined on its spectrum.
\endproclaim
\endremark

\head 5 Intrinsic Characterizations \endhead

In this section we use
 Proposition~17 to obtain an intrinsic
characterization result, Theorem~18; however,
the Proposition itself is perhaps of some interest.
Let $\sP\sI(\sD)$ denote the partial isometries in $\sD$.
We need the following lemma, which is equivalent to Lemma~3.5~(c) of
\cite{PePW1}; for completeness, we give a proof.

\proclaim{Lemma 16}
Let $\sA$ be a canonical algebra containing a canonical
masa $\sD$.
If $x$ and $y$ are elements of $N_\sD(\sA)$ with $\| x - y \| < 1$,
then $x^* y \in \sP\sI(\sD)$.
\endproclaim

\demo{Proof}
First, we show $x$ and $y$ have the same initial and final projections.
Suppose $x^* x$ and $y^* y$ are different.
Then there is a projection $p$ that is a subprojection of one and
orthogonal to the other.
Without loss of generality, assume $p \le x^* x$ and $p \perp y^* y$.
Then $xp = (x-y)p$ so $1 = \| xp \| = \| (x-y)p \| < 1$.
This proves that $x$ and $y$ have the same initial projections;
similarly, $x$ and $y$ have the same final projections.

Clearly, $x^*y$ is a partial isometry.
To prove $x^*y \in \sD$, it suffices to prove that $x^*y$ commutes with
all projections in $\sD$, or equivalently, with all subprojections of
$x x^* = y y^* $.
Let $p$ be such a projection.
We claim that $x p x^* = y p y^*$.
Accepting this for the moment, we have
$$ (x^*y) p = (x^* y) p (y^* y) = x^* (y p y^*) y = x^* (x p x^*) y =
(x^* x ) p (x^* y) = p (x^* y),$$
as required.

To prove the claim, suppose $x p x^* \ne y p y^*$.
Then there is some subprojection of one that is orthogonal to the other,
say $q$.
Without loss of generality, assume $q \le x p x^*$ and $ q \perp y p y^*$.
Then $qxp \ne 0$ and $qyp = 0$,
so $qxp = q(x-y)p$.
It follows that $1 = \| qxp \| = \| q (x-y) p \| < 1$;
a contradiction that proves the claim.
\qed
\enddemo
 \par
In the following Proposition, $\ra$ denotes the spectrum
for $\sA$ defined on the maximal ideal space of the
canonical masa $\sD$ and, for each
$e \in \nda$, $G(e)$ denotes the graph of the partial homeomorphism
induced by $e$ on the maximal ideal space of $\sD$.  Note that
$G(e)$ is a compact, open subspace of $\ra$.

\proclaim{Proposition 17}
Let $\sA$ be a canonical algebra containing a canonical
masa $\sD$.
Suppose $X \subset N_\sD(\sA)$ satisfies $\sP\sI(\sD) \cdot X \subseteq X$
$($or equivalently $ X \cdot \sP\sI(\sD) \subseteq X)$.
Then the following are equivalent:
\roster
\item The closed span of $X$ is $\sA$.
\item Each element of $N_\sD(\sA)$ can be written as a finite sum of
          elements of $X$.
\item In any presentation for $\sA$, each matrix unit can be written as
          a finite sum of matrix units in $X$.
\item $ \ra = \bigcup_{x \in X} G(x)$.
\endroster
If in addition $X \cdot X \subset X $, then we have another equivalent
condition:
\roster
\item"(5)" There is a presentation of $\sA$,  $\indlimit(B_i,\phi_i)$,
           with each matrix unit in each $B_i$ in $X$.
\endroster
\endproclaim

\demo{Proof}
$( 1 \Rightarrow 2)\quad$
Let $y \in N_\sD(\sA)$.
By hypothesis, the closed span of $X$ is $\sA$; so,
there is some $x$, a finite linear combination of elements
of $X$, such that $\| y - x \| < 1/2$.
Write $x$ as $\sum_{i=1}^l a_i x_i $, where $a_i$ is a scalar and
$x_i \in X$.
By rewriting the sum and restricting $x$ to the initial projection of $y$
(which we can do as $X \cdot \sP(\sD) \subseteq X$),
we may assume that each $x_i x_i^*$ is a subprojection of $y y^*$
orthogonal to all other $x_j x_j^*$'s.

If $\sum_{i=1}^l x_i x_i^*$ does not sum to $yy^*$, then letting
$z = yy^* - \sum_{i=1}^l x_i x_i^*$ we have $z(y- x) = zy$ and $\| z \| = 1$.
This gives a contradiction, since then
$\| y - x\| \ge \| z(y-x) \| = \| zy \| = 1$.
Thus $y y^* = \sum_{i=1}^l x_i x_i^*$, where each $x_i \in X$.
It follows that
$$
y = y y^* y = \left( \sum_{i=1}^n x_i x_i^* \right) y
= \sum_{i=1}^n x_i (x_i^* y).
$$

Since the initial and final projections of each $x_i$ are pairwise orthogonal,
$ \| x - y \| < 1/2 $ implies $ | 1 - a_i | < 1/2 $.
Letting $x' = \sum_{i=1}^l x_i$, it follows that
$ \| x' - y \| = \| x' -x \| + \| x - y \| < 1$.
As $x_i = x_i x_i^* x'$, we have
$$\| x_i - x_i x_i^* y \| \le \| x_i x_i^* \| \| x' - y \| = \| x' - y \| < 1$$
and so by Lemma~16, $x_i^* x_i x_i^* y = x_i^* y \in \sP\sI(\sD)$.
Since $x_i \in X$ and $x_i^* y \in \sP\sI(\sD)$, $x_i (x_i^* y) \in X$.
Thus $y$ is a finite sum of elements in $X$, as claimed.

\smallskip
\noindent
$( 2 \Rightarrow 3)\quad$
Let $\indlimit(A_i, \tau_i)$ be a presentation of $\sA$ and
let $t$ be a matrix unit in $A_a$ for some $a$.
By (2), $t=x_1 + x_2 + \cdots + x_n$ where each $x_n\in X$.
By Lemma~5.5 of \cite{Po4},
 every element of $N_\sD(\sA)$ (and in particular, each $x_n$)
can be written as a partial isometry in $\sD$ times
a finite sum of matrix units.
There is some $A_b$, $b \ge a$, that contains the finitely many matrix
units in the sums for $x_1,\ldots,x_n$.
Then we have $t = y_1 + y_2 + \cdots + y_m $ where each $y_j$ is
a partial isometry in $\sD$ times a matrix unit in $A_b$.
Since each $y_j$ is part of a sum that gives some $x_i \in X$,
by multiplying $x_i$ on the left by the projection $y_j y_j^*$
we have each $y_j \in X$.

We may assume that no sum of $y_a$'s equals zero by deleting all $y_a$'s
in such a sum.
Since $t$ also equals a sum of matrix units in $A_b$,
say $z_1 + \cdots + z_l$, we have that
$$ y_1 + \cdots + y_m - z_1 - z_2 - \cdots - z_l = 0. $$
This is only possible if $m=l$ and $y_1, \ldots, y_m$ is a permutation
of $z_1, \ldots, z_l$.
Thus $t$ is a finite sum of matrix units in $X$.
(Remark: here is another place where we use the assumption
that embeddings are regular.)

\smallskip
\noindent
$( 3 \Rightarrow 1)\quad$
This is immediate, as the closed span of $N_\sD(\sA)$ is $\sA$ and
(3) implies the closed span of $X$ contains $N_\sD(\sA)$.

\smallskip
\noindent
$(2 \Rightarrow 4)\quad$ Obvious from the definition of $\ra$.

\smallskip
\noindent
$(4 \Rightarrow 2)\quad$
Let $e \in \nda$.  Since $ G(e) \subset \ra
 = \bigcup_{x \in X}G(x)$ and $G(e)$ is compact, there exist
$x_1, \dots, x_n$ in $X$ such that $
G(e) \subseteq \bigcup_{i=1}^n G(x_i)$.  By multiplying the
$x_i$ by suitable projections from $\sD$, we may assume that
$G(e) = \bigcup_{i=1}^n G(x_i)$ and that the
$G(x_i)$ are pairwise disjoint.  It now follows that
$ e = \sum_{i=1}^n p_ix_i$ for some
$p_i \in \sP\sI(\sD)$.  By hypothesis, each $p_ix_i \in X$
and condition (2) holds.

\medskip
\noindent
We now also assume that $X \cdot X \subseteq X$ and prove that (5)
is equivalent to the first four conditions.

\smallskip
\noindent
$( 3 \Rightarrow 5)\quad$
Again, let $\indlimit(A_i, \tau_i)$ be a presentation of $\sA$.
For each $i$, let $X_i$ be the set of all matrix units in $A_i \cap X$.
Let $B_i$ be the closed span of $X_i$ for each $i$.
Observe that if $x,y \in X_i$ and $xy \ne 0$, then $xy \in X_i$
so $B_i$ is a subalgebra of $A_i$.

We claim that $\tau_i(B_i) \subset B_{i+1}$.
It suffices to show that if $x \in X_i$ and $y$ is a matrix
unit in $A_{i+1}$ that appears in the sum of matrix units $\tau_i(x)$,
then $y \in X_{i+1}$.
However, $x \in X_i$ implies $\tau_i(x) \in X$.
Since $X$ is closed under multiplication by projections in the
diagonal, the matrix unit $y =
  (y y^*) \tau_i(x) $ is in $X$ and so in $X_{i+1}$.

The hypothesis is that, for each $j$,  each matrix unit
 in $A_j$ can be written as a finite
sum of elements in $X$; hence $A_j \subset \cup_i B_i $.
This implies that
 $\overline{\cup_i B_i} = \overline{\cup_i A_i} = \sA$.
The presentation $\sA = \indlimit(B_i,\tau_i|_{B_i})$ now has
the required properties.

\smallskip
\noindent
$( 5 \Rightarrow 1)\quad$
We can repeat the (short) argument given in $(3 \Rightarrow 1)$.
\qed
\enddemo

Proposition~17 yields an intrinsic characterization of those TAF algebra
which have a presentation $\indlimit(A_i,\phi_i)$ where
each $\phi_{i,j} = \phi_j \circ \cdots \circ \phi_{i+1} \circ \phi_i$
is locally order preserving.

First, we need a few definitions.
Let $\sA$ be a subalgebra of an AF \cstar-algebra $\sB$ containing a
canonical masa $\sD$.
Following \cite{PePW1}, the diagonal ordering on the projections
in $\sD$, $\sP(\sD)$, is given by
$$ p \prec q \Longleftrightarrow \text{ there exists $w \in N_\sD(\sA)$
 with $ww^* = p$ and $ w^* w = q$}.$$

Given $w \in N_\sD(\sA)$, there is a partial homeomorphism given by
$x \mapsto w^* x w$ with domain $\{ x \in \sP(\sD) \stbar x \le w w^* \}$
and range $\{ x \in \sP(\sD) \stbar x \le w^* w \}$.
Call $w$ {\it order preserving\/} if this map preserves the diagonal
ordering restricted to its domain and range.
Let
$$ N^{op}_\sD(\sA) = \{ w \in N_\sD(\sA)
   \stbar w \text{ is order preserving} \}.$$
If $\sA=T_n$, this agrees with the definition of $N^{op}_{D_n}(T_n)$ in
Section~2.
Also note that $N^{op}_\sD(\sA)$ is intrinsic; it does
not depend on choosing
a presentation for $\sA$.

We should remark that $e \in N^{op}_\sD(\sA)$ if and only
$G(e) \times G(e)$, as a partial homeomorphism on $X \times X $,
sends $\ra$ into $\ra$.
(As usual, $X$ is the maximal ideal space of $\sD$.)
For subalgebras of groupoid \cstar-algebras, such subsets of
the spectrum (support subsemigroupoid) are called
{\it monotone $G$-sets} (page 57 of \cite{MS1}).
Groupoids admitting a cover of monotone $G$-sets (i.e., those
satisfying condition (3) below) have arisen in the study
of prime ideals; see Theorem~4.5 of \cite{MS1}.

\proclaim{Theorem 18}
Let $\sA$ be a canonical algebra containing a canonical
masa $\sD$.
The following are equivalent:
\roster
\item  The closed span of $N^{op}_\sD(\sA)$ is $\sA$.
\item  $\sA$ has a presentation $\indlimit(A_i,\phi_i)$
 so that, for all $i$ and $j$,
each $\phi_{i,j} = \phi_{j-1} \circ \cdots \circ \phi_i$
 is locally order preserving.
\item  $\ra =\bigcup \{G(e)\stbar e \in N^{op}_{\sD}(\sA)\}$.
\endroster
\endproclaim

\demo{Proof}
Since $N^{op}_\sD(\sA)$ satisfies all
 the conditions on $X$ in Proposition~17,
if the closed span of $N^{op}_\sD(\sA)$ is $\sA$ then
$\sA$ has a presentation $\indlimit(A_i,\phi_i)$ with each matrix
unit in $N^{op}_\sD(\sA)$.
It follows that each $\phi_{i,j}$ is locally order preserving.

The other direction follows from the observation that if all the
$\phi_{i,j}$ are locally order preserving,
then every matrix unit is in $N^{op}_\sD(\sA)$.

The equivalence of the condition that the spectrum
is the union of the graphs of the normalizing partial isometries
with the other two conditions follows immediately from
Proposition~17.
\qed
\enddemo

\remark{Remark} If $\sA$ is a strongly maximal TAF algebra, then
we can choose the presentation in condition~(2) to be strongly maximal,
i.e., each $A_i$ maximal triangular in $\roman{C}^*(A_i)$.
To see this, notice that in $(3 \Rightarrow 5)$ of the proof of
Proposition~17, if $A_i$ is maximal triangular in $\roman{C}^*(A_i)$ then
$B_i$ is maximal triangular in $\roman{C}^*(B_i)$.
Thus in Proposition~17, if $\sA$ is strongly maximal then it follows
that the presentation in condition~(5) can be chosen to be strongly maximal
and similarly in Theorem~18.
\endremark

The next two examples show that the class of algebras in Theorem~18
is properly contained in the class of algebras with
locally order preserving presentations and
properly contains the class of all algebras with
order preserving presentations using direct sums of $T_n$'s.
The first example appears in \cite{Do} as Example 13.
A similar example can be found in \cite{PePW2} (Example 3.7).

\example{Example}
Since locally order preserving embeddings are determined by
their action on the diagonal, we can define a locally
order preserving embedding $\phi_n \colon T_{3^n} \rightarrow
T_{3^{n+1}}$ by specifying the values of $\phi_n$ on
the minimal diagonal projections in $T_{3^n}$:
$$
\align
\phi_n(e_1^{(n)}) &= e_1^{(n+1)} + e_2^{(n+1)}
    + e_4^{(n+1)},\\
\phi_n(e_i^{(n)}) &= e_{3i-3}^{(n+1)} + e_{3i-1}^{(n+1)}
     + e_{3i+1}^{(n+1)}, \qquad \text{for } 1 < i < 3^n,\\
\phi_n(e_{3^n}^{(n)}) &= e_{3^{n+1}-3}^{(n+1)}
    +e_{3^{n+1}-1}^{(n+1)} + e_{3^{n+1}}^{(n+1)}.
\endalign
$$
Routine calculations will show that the composition of
two successive embeddings in this system fails to be
locally order preserving.  This alone is not sufficient
to show that the limit algebra obtained from this system
fails to satisfy the conditions of Theorem~18, since there
could, in  principle, be other presentations which satisfy
property (2).
\par
A matrix unit $e_{ij}^{(n)}$ in $T_{3^n}$ can be identified with
 its image in the limit algebra, $\sA$.  Observe that
$e_{ij}^{(n)} \in N^{op}_{\sD}(\sA)$ if, and only if,
$\phi_{m,n}(e_{ij}^{(n)})$ is order preserving for all $m>n$.
Using this, it is not hard to determine which matrix units
in $T_{3^n}$ are in the order preserving normalizer of $\sA$.
Indeed, it turns out that for $j>1$, $e_{1j}^{(n)}$ is not
in $N^{op}_{\sD}(\sA)$ and for $i<3^n$, $e_{i3^n}^{(n)}$
is not in $N^{op}_{\sD}(\sA)$.  All other matrix units
are in the order preserving normalizer.
\par
 From these observations, it is clear that the order
preserving normalizer fails to span $\sA$, so the algebra
$\sA$ is not in the family characterized in Theorem~18.
It is also illuminating to note that
$\bigcup \{G(e)\stbar e \in N^{op}_{\sD}(\sA)\}$
is an open, dense, proper subset of $\ra$.
\endexample

The next  example  shows that
$\overline{ \roman{span}\, N^{op}_\sD(\sA) } = \sA $
is not sufficient to imply the existence of a presentation with
order preserving embeddings between direct sums of $T_n$'s.
It has appeared before in the literature, as Example~3.2 in \cite{SV},
where it is shown to be a strongly maximal TAF algebra that is not analytic.

\example{Example}
Define $\phi_n \colon T_{2^n} \rightarrow T_{2^{n+1}}$ by
$$ \bmatrix A & B \\ & C \endbmatrix \mapsto
   \bmatrix A & & B & \\ & A & & B \\ & & C & \\ & & & C \endbmatrix
$$
where $A,C \in T_{2^{n-1}}$ and $B \in M_{2^{n-1}}$.
While $\phi_n$ is not order preserving
(consider $e_{1,1}+e_{2,1+2^{n-1}}$),
it does map an order preserving sum of matrix units in $A$ to an
order preserving sum.
Similar statements are true for order preserving sums in $B$ or $C$.

If $\sA = \indlimit(T_{2^n},\phi_n)$, then it is elementary to
see that each matrix unit is in $N^{op}_\sD(\sA)$ and so the closed span
of $N^{op}_\sD(\sA)$ is $\sA$.
As we noted above, \cite{SV} shows that this algebra is not analytic
and hence by Theorem~15, it cannot have a presentation using
order preserving embeddings through direct sums of $T_n$'s.
\endexample

\head 6 Intertwining Diagrams \endhead

Recall that a Banach algebra $\sA$ is the
{\sl inductive limit} of the system
$$
A_0 \symbarrow{\alpha_0} A_1 \symbarrow{\alpha_1} A_2 \symbarrow{\alpha_2}
A_3 \symbarrow{\alpha_3} A_4 \symbarrow{\alpha_4} \cdots
\tag{4}
$$
if there exists a sequence of injective homomorphisms
$\rho_i \colon A_i \rightarrow \sA$ so that
$ \sA = \overline{ \cup_i \rho_i(A_i) } $ and the diagram
$$ \CD
A_i @>\alpha_i>> A_{i+1} \\
& \symbse{\rho_i} & \symbdown{\rho_{i+1}} \\
& & \sA
\endCD  $$
commutes for every $i$.
Note that the injections $\rho_i$ and subalgebras $\rho_i(A_i)$ are
{\it not} unique in general.
For example, if $\alpha $ is an automorphism of $\sA$, then
each $\rho_i$ can be replaced by $\alpha \circ \rho_i$.
We will use this freedom in proving Theorem~19, below.

Also, note that if $\{n_i\}$ is any increasing sequence of
positive integers and $\alpha_{x,y} = \alpha_{y-1} \circ \cdots
\circ \alpha_{x+1} \circ \alpha_x$ then both (4) and
$$
A_{n_1} @>{\alpha_{n_1,n_2}}>> A_{n_2} @>{\alpha_{n_2,n_3}}>>
A_{n_3} @>{\alpha_{n_3,n_4}}>> A_{n_4} @>{\alpha_{n_4,n_5}}>> \cdots
 $$
have the same inductive limit.

The following theorem appears elsewhere in the literature
(Theorem~4.6 in \cite{V1}, Corollary~1.14 in \cite{PeW} and
the proof of Proposition~4.1 in \cite{DPo}) in a somewhat different form.
We include a proof for the sake of completeness; it is similar
to the proof in \cite{V1} but does not use the spectrum or cocycles.
Also, the proof has a straightforward extension to bimodules over
canonical masas.
Note that Theorem~19 is not applicable without the assumptions
of regularity and $*$-extendibility for the embeddings in the
direct limits.

\proclaim{Theorem 19}
Let $\sA = \indlimit(A_i, \alpha_i)$ and
 $\sB = \indlimit(B_i, \beta_i)$ be norm-closed
subalgebras of AF \cstar-algebras containing
 canonical masas $\sC$ and $\sD$
respectively.

If $\Phi \colon \sA \rightarrow \sB$ is an isometric algebra isomorphism
with $\Phi(\sC) = \sD$, then there exist strictly increasing sequences of
integers, $\{ m_i \}$ and $\{ n_i \}$ and
regular, isometric homomorphisms $\phi_i$ so that the diagram
$$ \CD
A_1 @>{\alpha_{1,m_1}}>> A_{m_1} @>{\alpha_{m_1,m_2}}>> A_{m_2}
      @>{\alpha_{m_2,m_3}}>> A_{m_3} @>{\alpha_{m_3,m_4}}>>
      A_{m_4} @>{\alpha_{m_4,m_5}}>> \cdots & \sA \\
& \symbse{\phi_1} & \symbup{\phi_2} & \symbse{\phi_3} & \symbup{\phi_4}
      & \symbse{\phi_5} & \symbup{\phi_6} & \symbse{\phi_7} & \symbup{\phi_8}
      & \symbse{\phi_9} & & \symbupdown{{\Phi,\Phi^{-1}}}\\
B_1 @>{\beta_{1,n_1}}>> B_{n_1} @>{\beta_{n_1,n_2}}>> B_{n_2}
      @>{\beta_{n_2,n_3}}>> B_{n_3} @>{\beta_{n_3,n_4}}>> B_{n_4}
      @>{\beta_{n_4,n_5}}>> \cdots & \sB
\endCD \tag{5} $$
commutes.
Moreover, however we identify each $A_i$ with an isomorphic subalgebra of $\sA$
we can identify each $B_i$ with an isomorphic subalgebra of $\sB$ so that
$\phi_{2i+1} = \Phi |_{A_{m_i}} $ and $\phi_{2i} = \Phi^{-1} |_{B_{n_i}}$.
\endproclaim

Given a diagram such as (5) with each $\phi_i$ isometric,
the universal property of inductive limits allows one to construct an
isometric isomorphism between the inductive limits.
What is useful to us is the conclusion that every isometric isomorphism
between TAF algebras arises in this way.

Also, notice that if $\sA$ and $\sB$ are triangular, then any
isometric isomorphism $\Phi \: \sA \rightarrow \sB$ necessarily
satisfies $\Phi( \sA \cap \sA^* ) = \sB \cap \sB^*$ and so
the assumption that $\Phi(\sC) = \sD$ is automatically satisfied.

\demo{Proof}
By the definition of inductive limit, there are nested subalgebras of $\sA$
and $\sB$ isomorphic to the algebras $A_i$ and $B_i$, respectively.
It is convenient to identify each $A_i$ and $B_i$ with its isomorphic
subalgebra.

Let $E$ be the set of matrix units in $\sA$ given by this identification and
similarly let $G$ be the set of matrix units in $\sB$.
Let $F = \Phi(E)$, a second set of matrix units in $\sB$.

Since $\Phi(\sC) = \sD$, it follows from Proposition~7.2 in \cite{Po4}
that $\Phi(N_\sC(\sA)) = N_\sD(\sB)$.
For each $g \in G$, Lemma~5.5 in \cite{Po4} implies that $g$ is
a partial isometry in $\sD$ times a sum of elements of $F$.
Suppose that $g = \delta_g f_g$ with $\delta_g$ a
 partial isometry in $\sD$
and $f_g$ a   sum of elements in $F$.
Without loss of generality, we may assume
$ g g^* = \delta_g \delta_g^* = \delta_g^* \delta_g = f_g f_g^*$ and
$ g^* g = f_g^* f_g$.

Define a map $\Gamma \colon G \rightarrow \sB$ by $\Gamma(g) = f_g$.
Notice that if $g = \sum_{i=1}^n g_i$ and $g,g_1,\ldots,g_n \in G$,
then $\delta_g f_g = \sum_{i=1}^n \delta_{g_i} f_{g_i}$.
As $g \in G$, the matrix units $g_i$ are pairwise orthogonal and hence so
are the matrix units $f_{g_i}$ and the partial isometries $\delta_{g_i}$.
Thus
$$ \delta_g f_g = \sum_{i=1}^n \delta_{g_i} f_{g_i} =
   \sum_{i=1}^n \delta_{g_i} \sum_{i=1}^n f_{g_i}.
$$
Multiplying by $\delta_g^*$ on the left and $\sum_{i=1}^n f_{g_i}^*$ on
the right, we have
$$ \delta_g^* \delta_g f_g \sum_{i=1}^n f_{g_i}^* =
   \delta_g^* \left( \sum_{i=1}^n \delta_{g_i} \right)
   \left( \sum_{i=1}^n f_{g_i} \right) \sum_{i=1}^n f_{g_i}^* .
$$
As $\delta_g^* \delta_g f_g = f_g f_g^* f_g = f_g$ and similarly
$ \sum_{i=1}^n \delta_{g_i} \sum_{i=1}^n f_{g_i} \sum_{i=1}^n f_{g_i}^* =
    \sum_{i=1}^n \delta_{g_i} $,
it follows that $f_g \sum_{i=1}^n f^*_{g_i} \in \sD$.
Since $f_g$ and $\sum_{i=1}^n f_{g_i}$ are sums of matrix units
and have the same initial and final projections,
$ f_g \ne \sum_{i=1}^n f_{g_i}$ would contradict
$f_g \sum_{i=1}^n f^*_{g_i} \in \sD$.
Thus $ f_g = \sum_{i=1}^n f_{g_i}$ or equivalently
$\Gamma(g) = \sum_{i=1}^n \Gamma(g_i)$.
It follows that we can extend $\Gamma$ by linearity to $\cup_i B_i$
in a well-defined way.

{\bf Claim:} $\Gamma \: \cup_i B_i \rightarrow \cup_i \Phi(A_i)$
is an isometric automorphism and a $\sD$-bimodule map.
Accepting this for the moment, it follows that $\Gamma$ can be extended
to an isometric automorphism of $\sB$.
Replacing each $B_i$ with the isomorphic subalgebra $\Gamma(B_i)$ does
not change the presentation of $\sB$.
However, by the definition of $\Gamma$,
any matrix unit in $F = \Phi(E)$ is a sum of matrix units in $\Gamma(G)$.
Conversely, the image under $\Phi^{-1}$ of any matrix unit in $\Gamma(G)$
is a sum of matrix units in $E$.

So with respect to the systems of matrix units $E$ and $\Gamma(G)$,
$\Phi$ and $\Phi^{-1}$ map matrix units to sums of matrix units.
Since there are only finitely many such matrix units in $A_1$ and each
is mapped to a finite sum of matrix units in $\sB$,
there is some $n_1$ so that $\Phi(A_1) \subset B_{n_1}$.
Continuing in this way, we have the required sequences and can
obtain each $\phi_i$ as the restriction of $\Phi$ or $\Phi^{-1}$.

It remains only to prove the claim.
By construction, $\Gamma$ is linear.
Suppose $g = g_1 g_2$ where $g_1,g_2 \in G$.
Then $ \delta_g f_g = \delta_{g_1} f_{g_1} \delta_{g_2} f_{g_2}$,
so
$$ f_g = ( \delta_g^* \delta_{g_1} ) f_{g_1} f_{g_1}^* f_{g_1}
           \delta_{g_2} f_{g_2}
       = ( \delta_g^* \delta_{g_1} ) f_{g_1} \delta_{g_2} f_{g_1}^*
           f_{g_1} f_{g_2}
       = ( \delta_g^* \delta_{g_1} ) ( f_{g_1} \delta_{g_2} f_{g_1}^* )
           f_{g_1} f_{g_2}.
$$
Hence $f_g f_{g_2}^* f_{g_1}^* \in \sD$.
Again $f_g$ and $f_{g_1} f_{g_2}$ are sums of matrix units with the
same initial and final projections so arguing as before,
we have $\Gamma( g_1 g_2 ) = \Gamma(g_1) \Gamma(g_2)$.
Since every element of $\cup_i B_i$ is a linear combination of elements
of $G$, it follows that $\Gamma$ is multiplicative.

Observe that if $g \in G \cap \sD$, then $g = f_g$ and $\Gamma(g) = g$.
To see this, we first observe that $g \in \sD$ and $f_g = \delta_g^* g$
so $f_g \in \sD$.
As $g$ and $f_g$ are projections, so is $\delta_g$ and hence
$f_g = f_g f_g^* = \delta_g^* \delta_g = \delta_g$.
Thus, $g = \delta_g f_g = f_g f_g = f_g$.
It follows that $\Gamma( d_1 b d_2 ) = d_1 \Gamma(b) d_2$ for
$d_1,d_2 \in \sD$ and $b \in \cup_i B_i$, so $\Gamma$ is a $\sD$-bimodule map.

If $f$ is a matrix unit in $\cup_i \Phi(A_i)$, then it can be written as
the product of a partial isometry in $\sD$ and a sum of matrix units in $G$,
say $f = \epsilon \sum_i g_i$.
On the other hand, each $g_i = \delta_{g_i} f_{g_i}$.
Since $\Gamma$ is a $\sD$-bimodule map,
$$\Gamma\left( \epsilon \sum_i \delta_{g_i} g_i \right)
    = \epsilon \sum_i \delta_{g_i} \Gamma\left( g_i \right)
    = \epsilon \sum_i \delta_{g_i} f_{g_i} = f,$$
so $\Gamma$ is surjective.

Since each pair $g_i$ and $f_{g_i}$ have the same initial and
final projections, we have
$$ \left\| \sum_{i=1}^n a_i g_i \right\| =
\left\| \sum_{i=1}^n a_i f_{g_i} \right\|;$$
thus $\Gamma$ is isometric.
Since $\Gamma$ is also surjective, it is an automorphism, as desired.
\qed
\enddemo

Notice that we have constructed an automorphism of $\sB$ that fixes
$\sD$ pointwise.
By Lemma~3.4 of \cite{V1}, such an automorphism is approximately inner.

\subhead Applications to Classifications \endsubhead
We outline the application of Theorem~19 to classifying
direct limit algebras with particular classes of embeddings.
While the refinement, standard, alternation, and
twist classifications given in  this section are well-known,
 these proofs seem simpler,
in part because of the common framework.
In each case, there are two key parts:
\roster
\item the $\phi_i$'s in the intertwining diagram~(5) have a nice form, and
\item embeddings in the class have a unique factorization.
\endroster
In Section 8 we give a new classification theorem for algebras
with order preserving embeddings.  In the next four examples  all
the finite dimensional algebras are full upper triangular matrix
algebras.

\example{Refinement Embedding Limit Algebras}
Consider the family of all refinement embeddings $\rho_k$.
Here, $k$ denotes the multiplicity of the embedding while
the dimension of the domain algebra is unspecified.
Since $\rho_k \circ \rho_l = \rho_{kl} = \rho_l \circ \rho_k$,
each refinement embedding can be factored as a composition of
refinement embeddings of prime multiplicity, and this factorization
is unique up to order.

Given $\indlimit(A_i,\alpha_i)$ with each $\alpha_i$ a refinement embedding
and $A_1=\Bbb{C}$, we can compute a supernatural number by, for each prime $p$,
counting the number of factors of multiplicity $p$ in the factorization
of $\alpha_{1,j}$ and taking the supremum as $j$ goes to infinity.

If $\indlimit(A_i,\alpha_i)$ and
 $\indlimit(B_i,\beta_i)$ are two direct limits
of this form and they have the same
 supernatural numbers as computed above, then it is routine
to construct an intertwining diagram such as (5).
It follows that the algebras are isometrically isomorphic.

On the other hand, suppose the two algebras are isometrically isomorphic;
we will show they have the same supernatural numbers.
By Theorem~19, we have an intertwining diagram such as (5).
Note that if $\gamma \circ \delta$ is a refinement embedding,
then necessarily $\delta$ is a refinement embedding, and hence
each $\phi_i$ in the diagram is a refinement embedding.

Since
$$ \alpha_{1,m_i} = \phi_{1,2i} \quad \text{ and } \quad
   \beta_{n_1,n_j} = \phi_{2,2j-1},
$$
it follows that the supernatural number of $\sA$ is given
by counting the refinement embeddings of each prime multiplicity in
$\phi_1,\phi_2,\ldots$ and the supernatural number of $\sB$ is
given by counting the number of refinement embeddings of each prime
multiplicity in $\beta_{1,n_1},\phi_2,\phi_3,\ldots$ .
Since $\phi_1$ and $\beta_{1,n_1}$ are both refinement embeddings
from $\Bbb{C}$ to $B_{n_1}=T_x$ for some $x$, they have the same
factorization.
Hence the supernatural numbers of $\sA$ and $\sB$ agree, as required.

This condition is necessary and sufficient, and so
classifies this family of algebras.
\endexample

\example{Standard Embedding Limit Algebras}
This classification proceeds in
 exactly the same way as for refinement embeddings.
\endexample

\example{Alternation Limit Algebras}
Suppose $\indlimit(A_i,\alpha_i)$ and $\indlimit(B_i,\beta_i)$ are
two direct limits with $A_1=B_1=\Bbb{C}$ and
each $\alpha_i$ and $\beta_i$ an alternation embedding
(a composition of standard embeddings and refinement embeddings).

Since $\rho_k \circ \sigma_l = \sigma_l \circ \rho_k$,
we can factor each alternation embedding as a composition of
standard embeddings and refinement embeddings, each of prime multiplicity.
Up to ordering, this factorization is unique.

For $\indlimit(A_i,\alpha_i)$, we can compute two supernatural numbers.
For the first, fix a prime $p$ and count the number of standard embeddings
of multiplicity $p$ in the unique factorization of $\phi_{1,j}$ for each $j$,
then take the supremum over all $j$.
Repeat this for each prime.
For the second, repeat this process, only counting refinement embeddings
of multiplicity $p$ instead of standard embeddings.
We will also consider the product of these two supernatural numbers,
which corresponds to counting all embeddings of each prime multiplicity,
standard and refinement.

First, if $\indlimit(A_i,\alpha_i)$ and $\indlimit(B_i,\beta_i)$
have their pairs of supernatural numbers agree, each up to a finite factor,
and the products of the pair are identical,
then the algebras are isometrically isomorphic.
Again, the construction of a diagram in the form of (5)
is routine.

To show the converse, we need the following fact:
If $\alpha \circ \beta$ is an alternation embedding,
then $\beta$ is an alternation embedding.
As $\alpha\circ\beta$ is an order preserving embedding,
Lemma~1 implies that $\beta$ is order preserving and so,
by Theorem~5, is a direct sum of refinement embeddings.
If the summands do not all have the same multiplicity (i.e., $\beta$ is
not an alternation embedding), then the summands in $\alpha \circ \beta$
will not all have the same multiplicity, a contradiction.

Suppose $\indlimit(A_i,\alpha_i)$ and $\indlimit(B_i,\beta_i)$ are
isometrically isomorphic.
Invoking Theorem~19, we have an intertwining diagram of
the form of (5).
By the previous paragraph, each $\phi_i$ in the diagram is an
alternation embedding.
We can compute a pair of supernatural numbers by counting the
refinement embeddings of each prime multiplicity and the standard
embeddings of each prime multiplicity in the sequence
$\phi_1,\phi_2,\ldots$ .
Since
$$ \alpha_{1,m_i} = \phi_{1,2i} \quad \text{ and } \quad
   \beta_{n_1,n_j} = \phi_{2,2j-1},
$$
each of the supernatural numbers for $\indlimit(A_i,\alpha_i)$
and $\indlimit(B_i,\beta_i)$ can differ by only a finite factor
from the supernatural numbers given by the alternation maps $\phi_i$.

In particular, each of the supernatural numbers for $\indlimit(A_i,\alpha_i)$
can differ by only a finite factor from the corresponding number for
$\indlimit(B_i,\beta_i)$.
Also, since the product of the two supernatural numbers corresponds to
counting all embeddings of a given prime multiplicity, the argument of
the previous examples shows that the products must agree exactly.
Thus the sufficient condition is also necessary.
\endexample

\example{Twist Embedding Limit Algebras}
A {\it twist} embedding is an embedding of the form
$\Ad\,U \circ \rho$, where $\rho$ is a refinement embedding
and $U$ is the permutation unitary matrix
which interchanges the last two minimal projections in the
diagonal.  In other words, $U$ is the identity  matrix with
the last two columns interchanged. Limit algebras constructed
with these embeddings were first
studied in \cite{PePW1}; their classification was given
in \cite{HPo}.
Unlike the previous examples, the composition of two twist
embeddings is not a twist embedding; the natural class to
consider consists of embeddings which are compositions of
twist embeddings.  If $\alpha\:T_k \longrightarrow T_{nk}$
is a composition of twist embeddings, then it is, in particular,
a nest embedding, i.e. it maps invariant projections under
$T_k$ to invariant projections for $T_{nk}$
\par
Quite generally, if $\phi$ is any embedding from $M_k$ into
$M_{nk}$, then we may write $\phi$ in the form $\Ad\,U \circ
\rho$ for some permutation unitary $U$.  The choice of $U$
is not unique.  Indeed, $\Ad\,U \circ \rho = \Ad\,V \circ \rho$
if, and only if, $\Ad\,V^{-1}U \circ \rho = \rho$; this
happens exactly when $V^{-1}U$ is block diagonal with
$n \times n$ blocks all of which are equal.
Also note that $\Ad\,U \circ \rho$ is a nest embedding if,
and only if, $U$ is block diagonal with each block of size
$n \times n$.  If we multiply each block on the right by a
fixed $n \times n$ permutation unitary matrix, then the resultant
matrix induces the same nest embedding as $U$ does.  This allows
us to replace $U$ by a matrix in standard form: multiply each
block on the right by the inverse of the first block.
So, we say that a block diagonal permutation matrix $U$
(with uniform block size) is in {\it standard form} if
the first block is the identity matrix.  One other trivial
fact about nest embeddings should be noted: if
$\phi \circ \psi$ is a nest embedding, then $\psi$ must be
a nest embedding.
\par
Suppose, now, that $\indlimit(A_i,\alpha_i)$ and
$\indlimit(B_i,\beta_i)$ are two direct limit algebras
with each $\alpha_i$ and $\beta_i$ a twist embedding.
By Theorem~19, we have an intertwining diagram as in (7).
For each $i$, $\phi_{i+1} \circ \phi_i$ is a nest embedding;
consequently, each $\phi_i$ is a nest embedding.  We shall
show below that each $\phi_i$ is a composition of twist
embeddings and further, that each composition of twist
embeddings has only one factorization into twist embeddings.
 From this we can conclude that there is some $m \le m_1$
such that $A_{m+i} = B_{n_1 + i}$ and $\alpha_{m+i} =
\beta_{n_1 + i}$, for all $i$.  This necessary condition
for isomorphic isomorphism is clearly also sufficient.
\par
In order to verify the second of the two claims above,
suppose that $\Ad\,V \circ \rho$ is a composition of
twist embeddings, where $V$ is in standard form.
The general observation: $\rho \circ \Ad\,U =
(\Ad\,\rho(U)) \circ \rho$ can then be used to see
that $V = U_q \circ \dots \circ U_1$
 is a product of permutation unitaries
each of which is the image of
an identity matrix with the last two columns interchanged
under a refinement embedding of suitable multiplicity.
 The critical observation is that it is possible to read
off from the matrix $V$, which is uniquely determined
by the requirement that it be in standard form,
the multiplicities of the refinements which are applied
to the $U_i$.  This yields a unique factorization for the
original embedding $\Ad\,V \circ \rho$ as a  composition
of twist embeddings.
\par
To verify the other claim, assume that
 $\tau = \Ad\,V \circ \rho_p
= \nu \circ \mu$, where $\tau$ is a composition of twist
embeddings, $V$ is a permutation unitary in standard form,
$\rho_p$ is a refinement of multiplicity $p$ defined on
some $T_k$, and $\nu$ and $\mu$ are nest embeddings with
multiplicities $n$ and $m$ respectively.  We need to prove
that $\nu$ and $\mu$ are actually compositions of twist embeddings.
Since $\nu$ and $\mu$ are nest embeddings, there are
unique permutation unitaries $V_n$ and $V_m$ in standard
form so that $\nu = \Ad\,V_n \circ \rho_n$ and
$\mu = \Ad\,V_m \circ \rho_m$.  Here, $\rho_n$ and $\rho_m$ are
refinement embeddings of multiplicities $n$  and $m$ and, of
course, $p = nm$.
 From the uniqueness of standard form and the general observation
in the preceding paragraph, it is easy to see that
$V = V_n\rho_n(V_m)$.  In order for $V$ to have the form
which the standard permutation unitary associated with
a composition of twist embeddings must have, it is necessary
that both $V_n$ and $V_m$ also have these forms.  This means
that $\nu$ and $\mu$ are compositions of twists.

\comment
Consider the family of ``refinement with twist'' embeddings,
that is embeddings $\Ad\,U \circ \rho$ where $U$ is an identity
matrix with the last two columns interchanged.
It is easy to show that any composition of such embeddings can be
factored uniquely in terms of such embeddings.

Suppose both $\delta \circ \gamma$ and $\gamma \circ \beta$ are
refinement with twist embeddings, then  $\gamma$ is an identity map.
To see this, let $U$ by an identity matrix with the last two
columns interchanged.
Then $\Ad U \circ \delta \circ \gamma$ is a refinement embedding
and so $\gamma$ is a refinement embedding.
Since $\gamma\circ\beta$ is a refinement with twist embedding,
it follows that $\gamma$ is an identity map.

Suppose $\indlimit(A_i,\alpha_i)$ and $\indlimit(B_i,\beta_i)$ are
two direct limits with $A_1=B_1=\Bbb{C}$ and
each $\alpha_i$ and $\beta_i$ a refinement with twist embedding.
By Theorem~19, we have an intertwining diagram such as (5).
Notice that both $\phi_{i+1}\circ\phi_i$ and $\phi_i\circ\phi_{i-1}$
are compositions of refinement with twist embeddings.
Since both have twists in the last embedding, the argument of
the previous paragraph shows that $\phi_i$ is an identity map.
Hence every $\phi_i$, $i > 1$ is the identity.
This implies that there is some $n$ so that $A_i = B_{i+n}$ and
$\alpha_i = \beta_{n+i}$ for sufficiently large $i$.

Clearly this necessary condition is also sufficient.
\endcomment
\endexample

\example{Ordered Bratelli diagrams with multiplicity}
As a final example, consider algebras $\indlimit(A_i,\alpha_i)$
where each $\alpha_i$ is order preserving
and each $A_i$ is a direct sum of $T_n$'s.
These algebras can be described in terms of ordered Bratelli diagrams,
introduced by Power in \cite{Po5}.
We begin by recalling the definition of ordered Bratelli diagram
given in \cite{PW}.  The definitions of ordered diagram
and ordered diagram with multiplicity have been given in
Section~2.

\definition{Definition}
An {\it ordered Bratelli diagram\/} is a pair $(\sV,\sE)$, where
$$\sV = V_0 \cup V_1 \cup \cdots,$$
a disjoint union of finite sets with $V_0$ a singleton, and
$$\sE = \{ (E_n,r_n,s_n) \stbar n \ge 1 \}$$
where each $(E_n,r_n,s_n)$ is an ordered diagram
from $V_{n-1}$ to $V_n$.

By {\it ordered Bratelli diagram with multiplicity}, we mean an
ordered Bratelli diagram as above with each $(E_n,r_n,s_n)$ replaced
by $(E_n,r_n,s_n,f_n)$, an ordered diagram with multiplicity.
\enddefinition

Using Theorem~6, we can associate an ordered Bratelli diagram
with multiplicity to each unital triangular AF algebra
$\indlimit(A_i,\alpha_i)$ where each $\alpha_i$ is order preserving,
each $A_i$ is a direct sum of $T_n$'s, and $A_0 = \Bbb{C}$.
We describe such triangular AF algebras in terms of their
ordered Bratelli diagrams with multiplicity.
This extends, in a natural way, Theorem~3.7 of \cite{PW} where
standard $\Bbb{Z}$-analytic TAF algebras are classified
by their associated ordered Bratelli diagrams.

First, we define the analogue of composition for ordered diagrams
with multiplicity, following \cite{PW}.

\definition{Definition}
Given two ordered diagrams with multiplicity,
$(E_1,r_1,s_1,f_1)$ from $V_1$ to $V_2$ and
$(E_2,r_2,s_2,f_2)$ from $V_2$ to $V_3$,
their {\it contraction\/} is an ordered diagram with multiplicity
$(E,r,s,f)$ from $V_1$ to $V_3$ given by
$$ E = \left\{ (e_1,e_2) \in E_1 \times E_2 \stbar r_1(e_1) = s_2(e_2)
            \right\}, $$
and
$$ s(e_1,e_2) = s_1(e_1), \quad r(e_1,e_2) = r_2(e_2), \text{ and }
   f(e_1,e_2) = f_1(e_1) f_2(e_2). $$
Given two edges, $(e_1,e_2)$ and $(f_1,f_2)$, with
$r(e_1,e_2) = r(f_1,f_2)$ then $r_2(e_2) = r_2(f_2)$.
If $e_2 \ne f_2$, then order $(e_1,e_2)$ and $(f_1,f_2)$
as $e_2$ and $f_2$ are ordered; if $e_2 = f_2$, then
$r_1(e_1) = r_1(f_1)$ and we can order $(e_1,e_2)$
and $(f_1,f_2)$ as $e_1$ and $f_1$ are ordered.

We denote $E$ as $E_2 \circ E_1$.
\enddefinition

As the notation suggests, if $\phi_1$ and $\phi_2$ are embeddings
associated to $E_1$ and $E_2$, then
$\phi_2 \circ \phi_1$ is the embedding associated to $E_2 \circ E_1$.
It follows that $(E_3 \circ E_2 ) \circ E_1 = E_3 \circ ( E_2 \circ E_1)$,
although this is also trivial to show directly.

Again following \cite{HPS} and \cite{PW}, we have:

\definition{Definition}
Given two ordered Bratelli diagrams with multiplicity, $(\sV,\sE)$
and $(\sW,\sF)$, we say they are {\it order equivalent\/} if
there exist strictly increasing functions
$f,g \colon \nbb \rightarrow \nbb$ and ordered diagrams with multiplicity
$E'_n$ from $V_n$ to $W_{f(n)}$ and $F'_n$ from $W_n$ to $V_{g(n)}$ so
that
$$ F'_{f(n)} \circ E'_n \cong^{ord} E_{g(f(n))} \circ \cdots \circ E_{n+1}
\tag{6} $$
and
$$ E'_{g(n)} \circ F'_n \cong^{ord} F_{f(g(n))} \circ \cdots \circ F_{n+1}
\tag{7} $$
for all $n \in \nbb$.
\enddefinition

It is now routine to prove:

\proclaim{Theorem 20}
Suppose $\sA = \indlimit(A_i,\alpha_i)$ and $\sB = \indlimit(B_i,\beta_i)$
are unital triangular AF algebras with each $A_i$ and $B_i$ a direct sums of
$T_n$'s, $A_0 = B_0 = \Bbb{C}$,
and each $\alpha_i$ and $\beta_i$ order preserving.

There is an isometric isomorphism $\Phi \colon \sA \rightarrow \sB$ if,
and only if, the ordered Bratelli diagrams with multiplicity associated to
$\indlimit(A_i,\alpha_i)$ and $\indlimit(B_i,\beta_i)$ are order equivalent.
\endproclaim

\demo{Proof}
Suppose there is an isometric isomorphism $\Phi \colon \sA \rightarrow \sB$.
By Theorem~19, we have the following commuting diagram:
$$ \CD
A_1 @>{\alpha_{1,m_1}}>> A_{m_1} @>{\alpha_{m_1,m_2}}>> A_{m_2}
      @>{\alpha_{m_2,m_3}}>> A_{m_3} @>{\alpha_{m_3,m_4}}>>
      A_{m_4} @>{\alpha_{m_5,m_4}}>> \cdots & \sA \\
& \symbse{\phi_1} & \symbup{\phi_2} & \symbse{\phi_3} & \symbup{\phi_4}
      & \symbse{\phi_5} & \symbup{\phi_6} & \symbse{\phi_7} & \symbup{\phi_8}
      & \symbse{\phi_9} & & \symbupdown{{\Phi,\Phi^{-1}}}\\
B_1 @>{\beta_{1,n_1}}>> B_{n_1} @>{\beta_{n_1,n_2}}>> B_{n_2}
      @>{\beta_{n_2,n_3}}>> B_{n_3} @>{\beta_{n_3,n_4}}>> B_{n_4}
      @>{\beta_{n_4,n_5}}>> \cdots & \sB .
\endCD \tag{8} $$
By Lemma~1, each $\phi_i$ is order preserving, since $\phi_{i+1} \circ \phi_i$
equals either some $\alpha_{a,b}$ or some $\beta_{a,b}$.
By Theorem~6, we can associate an ordered diagram with multiplicity
to each $\phi_i$, say $P_i$.

Let $(\sV,\sE)$ be the ordered Bratelli diagram with multiplicity
associated to $\indlimit(A_i,\alpha_i)$ and $(\sW,\sF)$ be the one
associated to $\indlimit(B_i,\beta_i)$.
To show they are order equivalent,
we define $f,g \colon \nbb \rightarrow \nbb$ by
$f(k) = n_j$ where $j$ is the least integer with $k \le m_{j-1}$ and
by $g(k) = m_j$ where $j$ is the least integer with $k \le n_j$.
Define $E'_k$ to be $P_{2j+1} \circ X$ where $X$ is the ordered diagram
with multiplicity associated to $\alpha_{k,m_j}$ and $j$ is the least
integer with $k \le m_j$.
Similarly, $F'_k$ is $P_{2j} \circ X$ where $X$ is the order diagram
with multiplicity associated to $\beta_{k,n_j}$ and $j$ is the least
integer with $k \le n_j$.
Commutativity of the diagram implies that (6) and (7) hold.

Conversely, if the diagrams are order equivalent, we can construct
sequences $\{m_i\}$ and $\{n_j\}$ and embeddings $\phi_i$ so that
the diagram (8) commutes.
It follows immediately that there is an isometric isomorphism
between $\sA$ and $\sB$.

Choose $n_1$ to be $f(1)$, $m_1$ to be $g(n_1)$, and for $i > 1$
choose $n_i$ to be $f(m_{i-1})$ and $m_i$ to be $g(n_i)$.
Define $\phi_1$ to be the embedding associated to the ordered
diagram with multiplicity $E'_1$, and for $i > 1$ define
$\phi_{2i}$ to be the embedding associated to $F'_{n_i}$
and $\phi_{2i+1}$ to be the embedding associated to $E'_{m_i}$.
Now (6) and (7) imply that the diagram (8) commutes.
\qed
\enddemo

Unlike the previous examples, here we
have found no way to pick a canonical
representative from an equivalence class
 of isometrically isomorphic algebras.
The difficulty is that we do not have a unique factorization theorem for
order preserving embeddings between direct sums of $T_n$'s.
Thus, Theorem~20 is a variant of Theorem~19 rather than a true
classification theorem.
In the next section, we restrict to order preserving embeddings between
$T_n$'s and obtain a unique factorization theorem.
Such a theorem is crucial for the classification given in
the final section.
\endexample

\head 7 Unique Factorization for Order Preserving Embeddings \endhead

In the last section, our aim was to demonstrate that
factorization theorems for families of embeddings yield
necessary and sufficient conditions for classification.
In this section, we prove unique factorization theorems for
order preserving embeddings.
While the proof is somewhat technical and requires several
preliminary lemmas, the statement of the factorization, Theorem~27, is simple.
In the next section, we will use this factorization to classify
limit algebras with order preserving presentations.

By Theorem~5, an embedding between $T_n$'s is order preserving
if, and only if, it is a direct sum of refinement embeddings.
For the sake of brevity, we use $\tupl a n$ to denote the direct sum of
$n$ refinement embeddings with multiplicities
$a_0, a_1, \ldots, a_{n-1}$ respectively,
that is, $ \rho_{a_0} \oplus \rho_{a_1}
 \oplus \cdots \oplus \rho_{a_{n-1}}$.
We refer to the number of entries in a tuple as its length,
and denote the length of $a$ by $\len a$.

With this notation, the composition of two embeddings, say
$a = \tupl a n $ and $ b = \tupl b m$, is an $mn$-tuple.
Notice that $ a \circ b $ equals
$$
( a_0 b_0, a_0 b_1, \ldots, a_0 b_{m-1}, a_1 b_0,
\ldots, a_1 b_{m-1}, a_2 b_0, \ldotss, a_{n-2} b_{m-1}, a_{n-1} b_0,
\ldots, a_{n-1} b_{m-1} )
$$
while in the other order, $ b \circ a $ equals
$$
( b_0 a_0, b_0 a_1, \ldots, b_0 a_{n-1}, b_1 a_0,
\ldots, b_1 a_{n-1}, b_2 a_0, \ldotss, b_{m-2} a_{n-1}, b_{m-1} a_0,
\ldots, b_{m-1} a_{n-1} ).
$$
Clearly, refinement embeddings commute with all embeddings.

Consider some tuple $a = \tupl a n$ with integer entries.
Dividing by $a_0$ gives a new, normalized tuple
$b = \tuplone b n$ with rational entries.  Because of
Theorem~13, we need only consider normalized tuples; i.e.,
tuples with rational entries and first entry always 1.

\pagebreak

One advantage of this standard form is that we immediately have
the following lemma:

\proclaim{Lemma 21}
If $b \circ c = \tuplone a n$,
and $\len c = m$, then $c = \tuplone a m $.

Hence if $b \circ c = d \circ e $ and $\len c = \len e$,
then $c = e$ and $b=d$.
\endproclaim

\demo{Proof}
That $c= \tuplone a m$ is immediate from the expression for composition.
It is obvious that $c = e$, while $b = d$ follows
from $c = e$ and the expression for composition.
\qed
\enddemo

Consider some tuple $a = \tuplone a n$. Given an integer
$m$, we say $a$ is {\it $m$-divisible\/} if $m$ divides $n$ and
the ratios $ a_i / a_{i-1} $ and $ a_j / a_{j-1} $ are equal
for all $i$ and $j$ such that $ i \equiv j \not\equiv 0 \pmod{m} $.
We will say $a$ is
{\it strongly $m$-divisible\/} if $m$ divides $n$ and
the ratios $ a_i / a_{i-1} $ and $ a_j / a_{j-1} $ are equal
for all $i$ and $j$ such that $ i \equiv j \pmod{m} $.
Notice that if $a = (1, x, x^2, x^3, \ldots x^{n-1})$ then $a$ will be
strongly $m$-divisible for $m$ any factor of $n$; in particular, $a$ is
strongly $1$-divisible if, and only if, $a$ is a geometric sequence.

\proclaim{Lemma 22}
Consider a tuple $a = \tuplone a n $  and
an integer $m$ such that $1 < m < n$.

Then there is a tuple $c$ such that $a = b \circ c$ with $\len c = m$
if, and only if, $a$ is $m$-divisible.
In this case, $c = \tuplone a m $ and
$b = (1, a_m, a_{2m}, a_{3m}, \ldots, a_{n-m} )$.

Further $a = b \circ c$ with $\len c = m$ and $b$ a geometric sequence
if, and only if, $a$ is strongly $m$-divisible.
In this case, $c = \tuplone a m $ and
$b = (1, a_m, a_{m}^2, a_{m}^3, \ldots, a_{m}^{n/m-1} )$.
\endproclaim

\demo{Proof}
Let $ k = n/m$.
We begin with $m$-divisibility.

If $a = b \circ c $ where $c$ is an $m$-tuple, then clearly $m$ divides $n$.
For any $l$ such that $0 \le l < k$, we have
$ ( a_{ml}, a_{ml+1}, \ldots, a_{ml+(m-1)} ) $ is $a_{ml}$ times $c$,
that is $ a_{ml} \tuplone a m$.
Thus $a$ is $m$-divisible.

If in addition, $b$ is a geometric sequence, then
$$ { a_{m(l+1)} \over a_{ml+(m-1)} } = { a_{m(l+1)} \over a_{ml}a_{m-1} }
= { a_m^{l+1} \over a_m^l a_{m-1} } = { a_m  \over  a_{m-1} }, $$
so $a$ is strongly $m$-divisible.

Conversely, suppose $a$ is $m$-divisible and let $b$ and $c$ be
as in the lemma.
If $i = qm+r$, with $0 \le r < m$, then the $i^{\roman{th}}$ entry
of $b \circ c$ is $b_q c_r$.
We will prove $b_q c_r = a_i$ by induction on $r$.
If $r=0$, then $ b_q c_0 = a_{qm} 1 = a_i $.
If the result holds for $r-1$, then $a_{i-1} = b_q c_{r-1} $.
Now as $ i \equiv r \not\equiv 0 \pmod{m} $, we have that
$ a_i /a_{i-1} = a_r / a_{r-1} $.
However by the definition of $c$, $ a_r /a_{r-1} = c_r / c_{r-1} $, so
$$ a_i = a_{i-1} { a_i \over a_{i-1} } = \left( b_q c_{r-1} \right)
{ a_r \over a_{r-1} } = b_q \left( c_{r-1} {c_r \over c_{r-1} }\right)
= b_q c_r. $$
Thus, $a = b \circ c$.

Strong $m$-divisibility implies that, for $j=1,2,\ldots,k-1$,
$$ { a_{jp} \over a_{(j-1)p + p-1} } = { a_p \over a_{p-1} } $$
and since $ a_{(j-1)p + p-1} = a_{(j-1)p} a_{p-1}$, we have
$a_{jp} = a_p a_{(j-1)p}$ for $j=1,\ldots,k-1$.
Thus $a_{jp} = a_p^j$, as required.
\qed
\enddemo

\proclaim{Lemma 23}
Let $p$ and $q$ be positive integers.
Suppose $a = \tuplone a n $ is $p$-divisible and $q$-divisible.
Then $a = b \circ c $ where $\len c = \lcm(p,q)$ and $c$ is
strongly $\gcd(p,q)$-divisible.
\endproclaim

\demo{Proof}
That $a$ can be factored as $b \circ c$
 with $\len c = \lcm(p,q)$ is immediate
from Lemma~22.
We need only show $c = (1,a_1,\ldots,a_{m-1})$ is
strongly $\gcd(p,q)$-divisible.

Suppose $i \in \{1,\ldots,p-1\}$ and $j \in \{1,\ldots,q-1\}$
such that $ i \equiv j \pmod{\gcd(p,q)}$.
It follows from the Chinese remainder theorem that
there is a unique integer $r$ less than $\lcm(p,q)$ so that
$r \equiv i \pmod{p}$ and $ r \equiv j \pmod{q}$.
The $p$-divisibility and $q$-divisibility imply that
$a_i/a_{i-1} = a_r/a_{r-1} = a_j /a_{j-1}$.
 From this, it is easy to show that if $i \equiv j \pmod{\gcd(p,q)}$
and either $i,j \in \{1,\ldots,q-1\}$ or $i,j \in \{1,\ldots,p-1\}$,
then $a_i/a_{i-1} = a_j /a_{j-1}$.

Turning to the general case, suppose $r$ and $s$ are integers so that
$0 < r,s < \lcm(p,q)$ and $r \equiv s \pmod{\gcd(p,q)}$.
There are integers $i_r$ and $j_r$ such that $0 \le i_r <p$
with $i_r \equiv r \pmod{p}$ and $0 \le j_r <q$ with
$j_r \equiv r \pmod{q}$.  As $0 < r < \lcm(p,q)$, at least
one of $i_r$ or $j_r$ is nonzero.  Similar statements
hold for $s$.
Thus, one of the four pairs
$$(i_r,j_s),\quad (j_r,i_s),\quad (i_r,i_s),\quad (j_r,j_s)$$
must have both entries nonzero.
Let $(k,l)$ be that pair; using $p$ and $q$-divisibility,
we have $a_r / a_{r-1} = a_k/a_{k-1}, a_s /a_{s-1} = a_l/a_{l-1}$,
and by the previous paragraph, $ a_k/a_{k-1} = a_l/a_{l-1}$.
Thus, $a_r/a_{r-1} = a_s/a_{s-1}$ and we are done.
\qed
\enddemo

\definition{Definition}
Call $a = \tuplone a n$ {\it irreducible} if $a$ cannot be factored
nontrivially.
\enddefinition

\proclaim{Lemma 24}
Given a tuple $a = \tuplone a n $, there is a
 minimal integer $m$ so that $a = b \circ c$ with $\len c = m$.
Further, $c = \tuplone a m $ and $c$ is irreducible.
\endproclaim

\demo{Proof}
Choose the least positive integer $m$, $1<m \le n$, so that
$a$ is $m$-divisible.
If $m = n$ then $a$ is
irreducible and we take $c = a$, $b = (1)$.
If $m<n$ then by Lemma~22
we can factor $a$ as $b \circ c$ with $\len c = m$.

If $c$ were reducible, this would contradict the minimality of $m$.
By Lemma~21, $c = \tuplone a m$ and we are done.
\qed
\enddemo

Geometric sequences do not have a unique factorization into irreducibles;
for example, $(1,x,x^2,x^3,x^4,x^5,x^6)$ can be written
as either $(1,x^3) \circ (1,x,x^2)$ or $(1,x^2,x^4) \circ (1,x)$.
In general, if $a$ is a geometric sequence and $a = b \circ c$,
then both $b$ and $c$ will be geometric sequences and the ratio of $b$
is exactly the ratio of $c$ raised to the power $\len c$.
It follows that for a tuple which is a geometric sequence, if we factor
its length into primes then for each distinct ordering of these primes
there will be a distinct factorization of the tuple.
There are two solutions to this problem: either we order such factors
according to length, or we avoid factoring such tuples at all.
We will give a factorization theorem for each of these solutions.

First we need a technical lemma.

\pagebreak

\proclaim{Lemma 25}
Suppose $x,y,z$ and $w$ are tuples with $y$ and $w$ irreducible
and $\len y \ne \len w$.

If $x \circ y = z \circ w$, then
$y$ is a geometric sequence and $x = z' \circ w'$ where $w'$ is
an irreducible geometric sequence so that $\len{w'}=\len{w}$
and $w' \circ y$ is a geometric sequence.
\endproclaim

\demo{Proof}
Let $m= \len{y}, n = \len{w}$, and $a = x \circ y$.
By hypothesis, $m \ne n$.

Notice that as $a$ is $m$-divisible and $n$-divisible, it is also
$\gcd(m,n)$-divisible.
Since $y$ and $w$ are initial segments of $a$
and $\gcd (m,n) < \max \{m,n\}$, if $\gcd(m,n) >1$ then
at least one of $y$ or $w$ is reducible, a contraction.

Thus $\gcd(m,n) = 1$ and so by Lemma~23 we may
conclude that $a = e \circ f$ where $\len f = m n $
and $f$ is a geometric sequence.
Also, $f = f_2 \circ f_1$ where $f_1$ and $f_2$ are geometric sequences,
$\len{f_2} = n$ and $\len {f_1} = m$.
Since $ e \circ f_2 \circ f_1  = x \circ y$ and $\len {f_1} = \len {y}$,
Lemma~21 shows $f_1 = y$ and $x = e \circ f_2$.

In particular, $y$ is a geometric sequence and if we let $z' = e$
and $w' = f_2$ then $w' \circ y$ is a geometric sequence.
Also, $\len{w'}  = n = \len{w}$, as required.
To show $w'$ is irreducible, we must first note that since $y$ is an
irreducible geometric sequence, $\len{y}$ must be a prime.
Since our hypothesis are symmetric in $w$ and $y$, similarly
$\len{w}$ is also a prime and so $w'$ is irreducible.
\qed
\enddemo

\proclaim{Theorem 26}
Every tuple has a unique factorization of the form
$c_1 \circ c_2 \circ  \cdots \circ c_k $, where each
$c_i$ is an irreducible tuple such that if
$c_i \circ c_{i+1}$ is a geometric sequence,
then $\len {c_i} \ge \len {c_{i+1}}$.
\endproclaim

\demo{Proof}
By repeatedly applying Lemma~24, we can factor $a$
into irreducibles and we claim that this factorization has the given form.
Suppose, for some $i$, $c_i$ and $c_{i+1}$ are geometric sequences
for which $c_i \circ c_{i+1}$ is also a geometric sequence.
Let $n = \len{c_i}$ and $m = \len{c_{i+1}}$.
We must show $n \ge m$.

Notice $c_i \circ c_{i+1} = e \circ f$ where
$f = (1,x,x^2,\ldots,x^{n-1})$ and
$e = (1, x^n, x^{2n}, \ldots, x^{(m-1)n})$.
Thus $ c_1 \circ c_2 \circ \ldots \circ c_{i+1}
 = c_1 \circ c_2 \circ \ldots c_{i-1} \circ e \circ f$.
However, when we factored $ c_1 \circ c_2 \circ \ldots \circ c_{i+1}$
using Lemma~24, we chose $c_{i+1}$ as
the factor of minimal length.
Hence $n = \len f \ge m$, as required.

\resultspace

It remains only to show that this factorization is unique.
Suppose $a$ can be factored in the above form
as $c_1 \circ c_2 \circ \cdots \circ c_k$
and as $d_1 \circ d_2 \circ \cdots \circ d_l$.
If $\len{c_k} = \len{d_l}$, then by Lemma~21, $c_k = d_l$ and
$c_1 \circ \cdots \circ c_{k-1} = d_1 \circ \cdots \circ d_{l-1}$.
By induction we are done, in this case.

Otherwise we have $\len{c_k} \ne \len{d_l}$ and will get a contradiction
by showing $\len{c_k} = \len{d_l}$.
By symmetry, it suffices to prove that $\len{d_l} \ge \len{c_k}$.

Since $\len{c_k} \ne \len{d_l}$, we can
apply Lemma~25 with
$$x = c_1 \circ \cdots \circ c_{k-1},\quad y = c_k,\quad
z = d_1 \circ \cdots \circ d_{l-1}, \hbox{ and } w=d_l.$$
Thus $c_k$ is a geometric sequence and
$c_1 \circ \cdots \circ c_{k-1} = z_1 \circ w_1$ with
$\len{w_1} = \len{d_l}$ and $w_1$ an irreducible geometric sequence
so that $w_1 \circ c_k $ is a geometric sequence.

If $\len{c_{k-1}} = \len{w_1}$, then by Lemma~21 $c_{k-1} = w_1$.
Thus $c_{k-1} \circ c_k$ is a geometric sequence.
By the assumed form of $c_1 \circ \cdots \circ c_k$, we can conclude
$\len{c_{k-1}} \ge \len{c_k}$ and so
$\len{d_l} = \len{c_{k-1}} \ge \len{c_k}$, as required.

Otherwise, $\len{c_{k-1}} \ne \len{w_1}$ so by Lemma~25 with
$$x = c_1 \circ \cdots \circ c_{k-2},\quad y = c_{k-1},\quad z = z_1,
\hbox{ and } w = w_1,$$
we have
$c_{k-1}$ is a geometric sequence and
$c_1 \circ \cdots \circ c_{k-2} = z_2 \circ w_2$, where
$\len{w_2} = \len{w_1}$ and $w_2$ an irreducible geometric sequence
so that $w_2 \circ c_{k-1} $ is a geometric sequence.

Further, since $w_1$ is a geometric sequence and
$c_1 \circ \cdots \circ c_{k-1} = z_1 \circ w_1$ we can conclude
that the ratio of $w_1$ equals the ratio of $c_{k-1}$.
As the ratio of $c_{k-1}$ equals the ratio of $w_2 \circ c_{k-1}$,
and $w_1 \circ c_k$ is a geometric sequence,
it follows that $w_2 \circ c_{k-1} \circ c_k $ is a geometric sequence.

If $\len{c_{k-2}} = \len{w_2}$, then $c_{k-2} = w_2$.
As before with $c_{k-1}$ and $w_1$, we may conclude that
$c_{k-2} \circ c_{k-1} \circ c_k $ is a geometric sequence and
$\len{d_l} = \len{c_{k-2}} \ge \len{c_{k-1}} \ge \len{c_k}$, as required.
If $\len{c_{k-2}} \ne \len{w_2}$, then we can apply Lemma~25 again.

Continuing in this way, we will either prove $\len{d_l} \ge \len{c_k}$
or end up with $c_1 = z_{k-1} \circ w_{k-1}$
where $\len{w_{k-1}} = \len{d_l}$.
Since $c_1$ is irreducible, it follows that $z_{k-1}= (1)$ and
so $w_{k-1} = c_1$.
As above, $c_1 \circ \cdots \circ c_k$ is a geometric sequence and
so $\len{c_1} \ge \len{c_2} \ge \cdots \ge \len{c_k}$.
Thus, $\len{d_l} \ge \len{c_k}$, as required.
\qed
\enddemo

It is clear that in the above factorization we can repeatedly multiply
together adjacent factors whose products will be geometric sequences
to get:

\proclaim{Theorem 27}
Every tuple has a unique factorization of the form
$c_1 \circ c_2 \circ \cdots \circ c_k $ satisfying:
\roster
\item"{a)}" each $c_i$ is either an irreducible tuple
             or a geometric sequence, and
\item"{b)}" for all $i<k$, $c_i \circ c_{i+1}$ is not a
             geometric sequence.
\endroster
\endproclaim

Finally, we consider when two embeddings commute:

\proclaim{Corollary 28}
Suppose $a$ and $b$ are normalized
tuples and $a \circ b = b \circ a$.
Then one of the following holds:
\roster
\item either $a= (1)$ or $b = (1)$,
\item $a = (1, \dots, 1)$ and $b = (1, \dots, 1)$, with possibly
       unequal lengths,
\item there exist a tuple $c$ and integers $m,n \in \Bbb{N} $
       so that $a=c^m$ and $b=c^n$.
\endroster
\endproclaim
\remark{Remark} In terms of commutativity of embeddings,
condition \therosteritem1 states that a refinement embedding
commutes with any order preserving embedding.  Condition
\therosteritem2 states that any two standard embeddings commute.
\endremark

\demo{Proof}
Let $d = a \circ b = b \circ a$.
If we are not in the third case, then it follows
that $d$ has two factorizations: one
given by factoring $a \circ b$ and another given by factoring $b \circ a$.
We can conclude by Theorem 26 that $d$, $a$ and $b$ are geometric series.

Further, $d = a \circ b$ implies the ratio of $a$
is the ratio of $b$ raised to the power $\len b$.
On the other hand $d = b \circ a$ implies the ratio of $b$
is the ratio of $a$ raised to the power $\len a$.
If $\len a = 1$ or $\len b = 1$, then we are in the first
case.  Otherwise, the ratios of $a$ and $b$ are both 1,
which is precisely the second case.
\qed
\enddemo

\head 8 Classification of Order Preserving Presentations \endhead

In this section,
we classify all  limit algebras which have
 presentations, $\indlimit(T_{n_i}, \alpha_i)$
with each $\alpha_i$ an order preserving embedding.
By Theorem~5, each $\alpha_i$ is a direct sum of
refinement embeddings.
Let
$\alpha_{x,y} = \alpha_{y-1} \circ \cdots \circ \alpha_{x+1} \circ \alpha_x$
for all integers $x$ and $y$ with $1 \le x < y$.
As in the last section, we can identify each $\alpha_{x,y}$ with a tuple,
the finite sequence of refinement multiplicities.

\definition{Definition}
We say that an order preserving presentation
of a limit algebra $\sA$ has
{\it geometric character\/} if
there is some $N$ so that for all $m$ and $n$ larger than $N$,
the tuple associated to $\alpha_{m,n}$ is a geometric sequence.
(We shall show that this is well-defined for
the limit algebra $\sA$ itself
in the course of proving the classification theorem.
Consequently, when $\sA$ has an order preserving presentation
with geometric character, we say that $\sA$ has geometric
character.)
\enddefinition

Geometric character implies that, for sufficiently large $m$,
the order preserving embedding
$\alpha_{1,m}$ factors as a finite sequence (not depending on $m$)
followed by a geometric sequence whose length depends on $m$.
Note that the ratio of this geometric sequence depends on the
choice of the finite sequence but not on $m$.

Choose the initial segment and consider the geometric sequences that follow it.
Since the product of an $m$-tuple and an $n$-tuple is an $mn$-tuple,
the length of the geometric sequence in $\alpha_{1,m}$ divides the length
of the geometric sequence in $\alpha_{1,m+1}$.
Thus, we can associate a supernatural number to this sequence of lengths
by counting the number of times a given prime divides any length in the
sequence.

By the {\it reduced root\/} of a rational number, $q$, we mean the
rational number $q^{1/n}$ such that for all $m > n$, $q^{1/m}$ is
not a rational number.
\par
Two supernatural numbers, $a$ and $b$, are {\it finitely
equivalent} if there are finite integers $m$ and $n$ so
that $ma$ and $nb$ are the same supernatural number.
In other words, if $a(p)$ and $b(p)$ are the exponents for the
prime $p$ in $a$ and  $b$ respectively, then
$a(p) \ne b(p)$ for only finitely many $p$ and only when
both $a(p)$ and $b(p)$ are both finite.
 Similarly, two unique
factorizations of sequences of normalized tuples are
{\it finitely equivalent} if either factorization can be
converted to the other by changing only finitely many factors.

If $\sA$ has geometric character, the invariants are:
\roster
\item the supernatural number of the \cstar-envelope (a UHF \cstar-algebra)
\item the finite equivalence class of the supernatural number of the first
     summands
\item the finite equivalence class of the supernatural number of the lengths
\item the reduced root of the ratio of the geometric sequence
\endroster

If $\sA$ does not have geometric character, the invariants are:
\roster
\item the supernatural number of the \cstar-envelope
     (again a UHF \cstar-algebra)
\item the finite equivalence class of the supernatural number of the first
     summands
\item the finite equivalence class of the unique factorization of the
     sequence of normalized tuples
\endroster
\remark{Remark}
Alternation algebras are a special subcase of the geometric
character case.  Each normalized tuple in the presentation for
an alternation algebra has all entries equal to 1.  The lengths
of the tuples are exactly the multiplicities of the standard
embedding factors.  Thus invariant (2) is just the finite
equivalence class of the refinement multiplicities, invariant
(3) is the finite equivalence class of the standard multiplicities,
and invariant (4) is necessarily equal to 1.
\endremark

\proclaim{Theorem 29}
Suppose $\sA$ and $\sB$ are triangular AF algebras and there is
an isometric isomorphism $\Phi \colon \sA \rightarrow \sB$.

If $\sA$ has a presentation, $\indlimit(T_{n_i},\alpha_i)$,
with each $\alpha_i$ order preserving,
 then so does $\sB$ and either they
both have geometric character or they both don't.
In either case, they have the same invariants.
\par
Conversely, two such algebras with the same invariants
are isometrically isomorphic.
\endproclaim

\demo{Proof}
Suppose that $\sA$
has a presentation $\indlimit(A_i,\alpha_i)$
with each $A_i$ the upper triangulars of some full matrix algebra and
each $\alpha_i$ a direct sum of refinement embeddings.
It is straightforward to see that $\sB$ has an order preserving
presentation,
since $\indlimit(A_i,\alpha_i)$ is also a presentation for $\sB$.
To see this formally, suppose for each $i$,
$\theta_i \colon A_i \rightarrow \sA$
is the isomorphism between $A_i$ and the isomorphic subalgebra of $\sA$.
The embedding $\Phi \circ \theta_i\: A_i \rightarrow \sB$
gives a subalgebra of $\sB$ isomorphic to $A_i$, for each $i$.

Suppose $\sB$ has some other presentation $\indlimit(B_i,\beta_i)$
with each $B_i$ the upper triangulars of some full matrix algebra and
each $\beta_i$ a direct sum of refinement embeddings.
We must show that this presentation of $\sB$ has geometric character
if and only if the presentation of $\sA$ does.

By Theorem 19, we have an intertwining diagram:
$$ \CD
A_1 @>{\alpha_{1,m_1}}>> A_{m_1} @>{\alpha_{m_1,m_2}}>> A_{m_2}
    @>{\alpha_{m_2,m_3}}>> A_{m_3} @>{\alpha_{m_3,m_4}}>>
    A_{m_4} @>{\alpha_{m_5,m_4}}>> \cdots & \sA \\
& \symbse{\phi_1} & \symbup{\phi_2} & \symbse{\phi_3} & \symbup{\phi_4}
      & \symbse{\phi_5} & \symbup{\phi_6} & \symbse{\phi_7} & \symbup{\phi_8}
      & \symbse{\phi_9} & & \symbupdown{{\Phi,\Phi^{-1}}}\\
& & B_{n_1} @>{\beta_{n_1,n_2}}>> B_{n_2} @>{\beta_{n_2,n_3}}>> B_{n_3}
        @>{\beta_{n_3,n_4}}>> B_{n_4} @>{\beta_{n_4,n_5}}>> \cdots & \sB
\endCD \tag{9}$$

If the presentation of $\sA$ has geometric character, we may choose
the sequence $\{m_i\}$ so that
$\alpha_{m_i,m_j}$ is a geometric sequence for each $i$ and $j$ with
$1 \le i < j$.
It is straightforward to observe that the product of two tuples
is a geometric sequence only if both of the original tuples are
geometric sequences.
Since $\alpha_{m_i,m_{j}} =
\phi_{2j} \circ \beta_{n_{i+1},n_{j}} \circ \phi_{2i+1}$,
it follows that $\beta_{j,k}$ is a geometric sequence for all $j$ and $k$
with $k > j \ge n_2$.
Thus, the presentation of $\sB$ has geometric character.

Similarly, if the presentation of $\sB$ has geometric character,
then so does the presentation of $\sA$.

We now prove that $\sA$ and $\sB$ have the same invariants.

Using Theorem~7.5 of \cite{Po4} we can extend $\Phi$ to a $*$-isomorphism
between the \cstar-envelopes of $\sA$ and $\sB$.
Since the \cstar-envelopes are UHF \cstar-algebras, by Glimm's classification
the supernatural numbers of the \cstar-envelopes must agree.

We have already proved, in Theorem 13, that the supernatural numbers of the
first refinement summands must agree up to finite equivalence.
We now divide all tuples through by their
first entry, so that all tuples
begin with $1$.
This allows us to apply the unique factorization theorem of the last section,
Theorem~27.

\medskip
\noindent
{\sl Case 1: $\sA$ has geometric character.}
 From the diagram (9), we have
$$ \alpha_{m_1,m_j} =
 \phi_{2j} \circ \beta_{n_2,n_j} \circ \phi_3. \tag{10}$$
As above, we may assume that $\alpha_{m_1,m_j}$ is a geometric sequence
for each $j > 1$.
It follows that $\beta_{n_2,n_j}$ is also a geometric sequence for each
$j \ge 2$.

Recall that for any tuples $a$ and $b$, the product $b \circ a$ is a geometric
sequence only if both $a$ and $b$ are geometric sequences with the ratio
of $b$ equal to the ratio of $a$ raised to the power $\len\,b$.
Thus the ratio of $\beta_{n_2,n_j}$ is a power of the ratio of
 $\alpha_{m_1,m_j}$, and therefore a power of the reduced root
of the presentation for $\sA$.

Observe that $\beta_{1,n_j} = \beta_{n_2,n_j} \circ \beta_{1,n_2}$ will
factor as some initial segment followed by a geometric sequence.
The ratio of this geometric sequence will be a rational number that is a
root of the ratio of $\beta_{n_2,n_j}$.
In particular, we can conclude that the reduced roots given by these two
presentations of $\sA$ and $\sB$ must be the same.

Also from (10), we can conclude that, for each $j>2$,
the length of the geometric sequence
$\beta_{n_2,n_j}$ must divide the length of the geometric sequence
 $\alpha_{m_1,m_j}$.
It follows that, after deletion of a finite factor,
 the supernatural number associated    with the
lengths of the geometric sequences for $\sB$ divides the
 corresponding supernatural
number for $\sA$.
Since we may interchange the roles of $\sA$ and $\sB$, we have
that, up to finite equivalence, the supernatural numbers associated
with the lengths agree.
This completes case~1.

\medskip
\noindent
{\sl Case 2: $\sA$ does not have geometric character.}
Since $\sA$ does not have geometric character, for any integer $y$ there is
an integer $z > y$ so that the (unique) factorization of $\alpha_{y,z}$ is not
a geometric sequence.
Hence for any $\alpha_{x,y}$, there is a $z > y$ so that for any $w > z$,
the factorization of $\alpha_{x,w}$ is
the factorization of $\alpha_{x,z}$ followed by
the factorization of $\alpha_{z,w}$.
As $\sB$ also does not have geometric character, it follows that
a similar statement holds for the maps $\beta_{x,y}$.

Consequently, we can choose the sequences
 $\{m_i\}$ and $\{n_i\}$ so that
the factorization of $\alpha_{m_j,m_k}$ is
the factorizations of $\alpha_{m_i,m_{i+1}}$ for $j \le i < k$ in order
and similarly for the $\beta_{n_j,n_k}$.
Since
$$\beta_{n_k,n_{k+1}} \circ \beta_{n_{k-1},n_k} =
      ( \phi_{2k+1} \circ \phi_{2k} ) \circ
 ( \phi_{2k-1} \circ \phi_{2k-2} ) $$
$$\alpha_{m_k,m_{k+1}} \circ \alpha_{m_{k-1},m_k} =
      ( \phi_{2k+2} \circ \phi_{2k+1} ) \circ
 ( \phi_{2k} \circ \phi_{2k-1} ),$$
it follows  that the unique factorization of $\phi_k \circ \phi_{k+1}$
is the unique factorization of $\phi_k$ followed by that of $\phi_{k+1}$
for every $k > 1$.
Since $\alpha_{m_k,m_{k+1}} = \phi_{2k+2} \circ \phi_{2k+1}$
and $\beta_{n_k,n_{k+1}} = \phi_{2k+1} \circ \phi_{2k}$
 for every $k > 1$,
it follows that (after removing all $1$-tuples) the
sequences of tuples
in the presentations of $\sA$ and $\sB$ have the
 same unique factorization, except for possibly
different initial segments ($\alpha_{1,m_1}$ and $\beta_{1,n_1}$).

This completes case~2 and the classification.
\qed
\enddemo

\enddocument